\documentclass[10pt]{amsart}
\usepackage{latexsym, amsmath, amsfonts, amssymb, mathrsfs, fancyhdr, tikz, tikz-cd}
\usepackage{verbatim}
\usepackage{multicol}
\usepackage{enumitem} 
\usepackage{appendix}
\usetikzlibrary{matrix,arrows,decorations.pathmorphing,calc}
\usepackage[normalem]{ulem}

\def\R{\mathbb{R}}

\def\v{\varphi}
\def\e{\varepsilon}
\def\s{\sigma}

\def\d{\delta}

\def\l{\lambda}

\def\S{\mathcal{S}}

\def\a{\alpha}
\def\b{\beta}

\def\t{\tau}

\def\N{\mathbb{N}}
\def\la{\langle}
\def\ra{\rangle}

\def\fe{f_\e}

\def\ch2{\mathbb{C} \mathbb{H}^2}

\def\th{\theta}

\def\oshr/2{\cosh \left( \frac{r}{2} \right)}
\def\inhr/2{\sinh \left( \frac{r}{2} \right)}
\def\osh2r/2{\cosh^2 \left( \frac{r}{2} \right)}
\def\inh2r/2{\sinh^2 \left( \frac{r}{2} \right)}
\def\-1/4{- \frac{1}{4}}
\def\H{\mathbb{H}}

\def\S{\mathbb{S}}

\def\Z{\mathbb{Z}}
\def\h{\mathfrak{h}}

\def\G{\Omega}
\def\fe{f_{\epsilon}}
\def\ve{\varphi_{\epsilon}}
\def\fd{f_\delta}

\DeclareMathOperator{\SL}{SL}
\DeclareMathOperator{\SO}{SO}
\DeclareMathOperator{\Isom}{Isom}

\DeclareMathOperator{\Vol}{Vol}

\newtheorem{theorem}{Theorem}[section]
\newtheorem{lemma}[theorem]{Lemma}
\newtheorem{cor}[theorem]{Corollary}

\newtheorem{proposition}[theorem]{Proposition}

\theoremstyle{definition}
\newtheorem{definition}[theorem]{Definition}
\newtheorem{example}[theorem]{Example}

\theoremstyle{remark}
\newtheorem{remark}[theorem]{Remark}

\numberwithin{equation}{section}

\begin{document}

\title{Codimension two complements in hyperbolic products}

\author{Hyeran Cho, Jean-Fran\c cois Lafont, Nicholas Miller, Barry Minemyer}

\date{\today.}

\begin{abstract}
Let $M$ be a finite-volume locally symmetric manifold modeled on $\H^m \times \H^n$, and let $S$ be a compact, totally geodesic, codimension two submanifold of $M$ modeled on $\H^{m-1} \times \H^{n-1}$. In this paper, we initiate the study of manifolds of the form $N = M \setminus S$. 
In the case where $M$ is a product and $S$ is connected, we prove that such $N$ admit complete, finite-volume, $A$-regular metrics with non-positive sectional curvature. We also deduce some topological properties of $N$, and algebraic properties of $\pi_1(N)$. In the case where $M$ is irreducible, we show that $N$ admits {\it almost non-positive curvature}. Similar results hold for manifolds $X_d$, which are $d$-fold branched covers of $M$ ramified over $S$.
Along the way, we show that there are infinitely many commensurability classes of such pairs $(M,S)$, and develop general curvature formulas for doubly warped product metrics, which may be of independent interest. 
\end{abstract}

\maketitle

\section{Introduction}\label{Section:Introduction}

In the study of Riemannian manifolds of non-positive curvature, a persistent issue is the relative paucity of examples. One common theme is to find procedures to create new examples from existing examples. To this end, assume $M$ is a finite-volume, complete, non-positively curved Riemannian manifold, and $S\subset M$ is an embedded totally geodesic, codimension two submanifold. Let $N = M \setminus S$ and, for some integer $d > 2$, suppose that there is a smooth manifold $X_d$ which is a $d$-fold ramified branched cover of $M$ about $S$.  

\vskip 10pt

\noindent{\bf Key Problem:} From a differential geometric point of view the central questions are the following:
\begin{enumerate}
\item Does $N$ admit a complete, finite-volume, non-positively curved metric?
\item Does $X_d$ admit a non-positively curved Riemannian metric?
\end{enumerate}

\vskip 10pt

These questions have been studied in detail when $(M, S)$ is modeled on a {\it negatively curved} symmetric space and a totally geodesic codimension two symmetric space. For the branched cover construction $X_d$, the first results were obtained by Gromov--Thurston \cite{GromovThurston} in the case where $M$ is real hyperbolic, with the pair modeled on $(\H^n, \H^{n-2})$. They constructed examples where the singular pull-back metric could be smoothed to have curvature arbitrarily close to $-1$. Recently, Stover--Toledo \cite{StoverToledo} constructed branched covers over some complex hyperbolic manifolds, and Minemyer \cite{MinemyerKahler} resolved the smoothing problem, putting negatively curved metrics with pinching constant close to $1/4$ on these branched covers. This corresponds to the case $(\mathbb C\H^n, \mathbb C \H^{n-1})$. Several other interesting Riemannian smoothings have been constructed on the Stover--Toledo examples, see e.g. Zheng \cite{Zheng}, Guenancia--Hamenst\"adt \cite{GH-KahlerEinstien}, and Lafont--Minemyer \cite{LafontMinemyer}.

For the complement $N$, the real hyperbolic case was studied by Belegradek \cite{belegradek2007rigidity}, while the complex hyperbolic case modeled on $(\mathbb C \H^n, \mathbb C \H^{n-1})$ was worked out in a subsequent paper \cite{belegradek2012complex}. In both these cases, complete finite-volume negatively curved metrics exist on the complement $N$. The exceptional case modeled on the totally real embedding $(\mathbb C \H^2, \H^2)$ was studied by Minemyer \cite{minemyer2018real}, who constructed the desired metrics.

Moving to the case of non-positive curvature, the questions become much more challenging. The key issue is that negative curvature is an open condition, while small perturbations of a non-positively curved metric often introduces positive curvature. A breakthrough result by Stadler \cite{Stadler} shows that for the branched cover case modeled on the diagonal embedding $(\H^2 \times \H^2, \Delta(\H^2))$, examples exist where the branched cover $X_d$ does {\bf not} support any non-positively curved Riemannian metric. This is the only case where this question has been considered.

\vskip 10pt

\subsection*{General case}

In the present paper, we initiate the study of the problem in the case where the pair is modeled on $(\H^m\times \H^n, \H^{m-1} \times \H^{n-1})$ for $m,n\ge 2$. While we are unable to resolve the original problem, we can show that the smooth manifolds $X_d$ and $N$ are {\it almost non-positively curved}. 
Our main result is the following:

\begin{theorem}\label{thm:main theorem 2}
Let $(M, S)$ be locally modeled on the pair $(\H^m\times \H^n, \H^{m-1} \times \H^{n-1})$, and let $N$ and $X_d$ be as defined above. Then for any $\e > 0$ we have
\begin{enumerate}
\item the non-compact manifold $N$ admits a complete Riemannian metric $g_\e$ with the property that all sectional curvatures are $\leq \e$, and $\Vol(g_\e)\leq V$ for some uniform constant $V$ (independent of $\e$). 
\item the branched cover $X_d$ admits a Riemannian metric $h_\e$ with the property that all sectional curvatures are $\leq \e$. Moreover, given any $\xi>0$, the metric $h_\e$ can be chosen to satisfy $\Vol(X_d, h_\e)\leq d\Vol(M) + \xi$.
\end{enumerate}
\end{theorem}

Note that in the theorem above, the uniform upper bound $V$ on volume is used as a normalizing condition. Otherwise one could just take any complete Riemannian metric on $M$, and scale the metric in order to make the sectional curvature $\leq \epsilon$.

In both of these constructions, the approach consists of modifying the symmetric metric in a very small neighborhood of $S$, using a warped product construction. The non-positive curvature constraint then turns into a system of non-linear partial differential inequalities. Exploiting the symmetries of the construction, this can be reformulated as a condition on a certain family of two-by-two matrices to be negative semi-definite. Analyzing the spectral properties of this family of matrices leads us to our result.

These smoothing problems are quite delicate. For an illustrative example, consider the hyperelliptic involution, realized as an isometric involution on a flat $2$-torus. The quotient space is topologically a $2$-sphere, equipped with a flat metric with four singular points of cone angle $\pi$. While it is possible to smooth these four singularities in an arbitrarily small neighborhood of the singular points, Gauss--Bonnet then forces the curvature in that neighborhood to tend to $+\infty$. Indeed, Gauss--Bonnet implies that $S^2$ is 
{\bf not} almost non-positively curved.

In contrast, our result tells you that, if you allow yourself an arbitrarily small amount of positive curvature, you can smooth the singularity in the branched cover $X_d$. Our smoothing can also be carried out in an arbitrarily small neighborhood of the branching locus $S$ (though in this case we will have tangent $2$-planes where the sectional curvature tends to $-\infty$). The interested reader can also compare this to the Gromov--Thurston examples \cite{GromovThurston} and the Stover--Toledo examples \cite{MinemyerKahler}, where smoothing while keeping curvature pinched close to that of the symmetric metric requires a large normal injectivity radius around the singular locus $S$. 

\vskip 10pt

\subsection*{Product case}
Since we were unable to construct a non-positively curved smoothing working purely locally, it is natural to consider adding some global hypotheses. 
 
When $M$ is a product, the question becomes much more tractable, and we can show:

\begin{theorem}\label{thm:main theorem 1}
Let $(M, S)$ be locally modeled on the pair $(\H^m\times \H^n, \H^{m-1} \times \H^{n-1})$, and assume that $M$ is a product, and $S$ is connected. If $N$ and $X_d$ are defined as above, then
\begin{enumerate}
\item $N$ admits a complete, finite-volume, $A$-regular\footnote{See Definition \ref{def:A-regular}.}, non-positively curved Riemannian metric, and
\item $X_d$ admits a non-positively curved Riemannian metric.  
\end{enumerate}
\end{theorem}

Statement (2) above is due to Fornari--Schroeder \cite{fornari1990ramified}, and is included for completeness. The idea here is to leverage the global product structure to construct a non-positively curved metric on $M$ that is locally flat transverse to the subset $S$, i.e. so that transverse to $S$ the metric is locally isometric to $\H^{m-1} \times \H^{n-1}\times \R ^2$. Fornari--Schroeder then show that the pull-back metric on $X_d$ can be locally smoothed near $S$ while keeping non-positive curvature. For statement (1), we adapt the $A$-regular metric construction of Belegradek in \cite{belegradek2012complex} to show how to modify the transversely locally flat metric in a (punctured) tubular neighborhood of $S$ to create the desired cusp. This is similar to the ``funnel" construction of Avramidi and Nguyen-Phan in \cite{AvramidiPham}. 
Theorem \ref{thm:main theorem 1} holds in greater generality.  The submanifold $S$ does not need to be connected as long as its projections onto the factors are disjoint.

\vskip 10pt

\subsection*{Topological and algebraic applications}

When a manifold supports metrics of non-positive curvature, there are some strong topological and algebraic consequences. As a first application, observe that the manifold $N$ in Theorem \ref{thm:main theorem 1} can be compactified to an aspherical manifold with boundary $\partial\bar N$. Combining Theorem \ref{thm:main theorem 1} with work of Farrell--Jones \cite{FarrellJones} allows us to conclude a topological rigidity result:

\begin{cor}[Borel Conjecture]\label{cor:rigidity of N}
    Suppose $(M, S)$ is as in Theorem \ref{thm:main theorem 1}, and assume $\dim M \geq 5$. Let $N'$ be an  aspherical manifold with boundary $\partial N'$, and assume $\phi: (N', \partial N') \rightarrow (\bar N, \partial \bar N)$ is a homotopy equivalence which is a homeomorphism on the boundary. Then $\phi$ can be homotoped rel boundary to a homeomorphism of pairs.
 \end{cor}

From the group theoretic viewpoint, the relationship between the fundamental group of punctured manifolds $N$ and the fundamental group of $M$ is complicated. The $\pi_1(N)$ are examples of {\it overlattices} in the sense of Gromov \cite[Pg 126]{Gromov}, groups that come with a natural surjection onto a uniform lattice. In the case where $M$ is negatively curved, and $S$ is totally geodesic, Belegradek--Hruska \cite{BelegradekHruska} showed $\pi_1(N)$ is relatively hyperbolic, allowing them to deduce a large number of properties for $\pi_1(N)$ (see also \cite{belegradek2007rigidity}, \cite{belegradek2012complex}). In contrast, in our setting the pair $(\pi_1(\bar N), \pi_1(\partial \bar N))$ cannot be relatively hyperbolic due to the abundance of $\Z^2$-subgroups away from the ends. Nevertheless, in the product case we have certain graph of groups decompositions, which we can leverage to verify some algebraic properties for $\pi_1(N)$.

\begin{theorem}[Algebraic Properties]\label{thm: group theoretic properties}
Assume $M$ is compact, a product, and that $S$ is connected. Then for the complement $N$, $\pi_1(N)$ satisfies the following properties:
\begin{enumerate}
\item uniform exponential growth
\item Tits alternative
\item Kazhdan subgroups are always trivial
\item finite asymptotic dimension
\item Baum--Connes conjecture (with coefficients)
\item property RD (rapid decay)
\end{enumerate}
\end{theorem}

\subsection*{Layout of the paper}
The remainder of this paper is laid out as follows.  
In Section \ref{sec:background}, we give preliminaries on locally symmetric spaces modeled on $\H^m\times \H^n$, culminating in Theorem \ref{thm:construction}, which shows that for $m, n\geq 3$, there are infinitely many commensurability classes of manifolds fitting the hypothesis of Theorem \ref{thm:main theorem 2} .
In Section \ref{Section:Prelim}, we explain and add details to a construction that was quickly outlined \cite[Sec 3]{fornari1990ramified} and use this to prove Theorem \ref{thm:main theorem 1} (1).  
Additionally, we give more details as to why this argument fails when $M$ is not a product.
We also develop curvature formulas for the warped product metric corresponding to the product metric $\tilde{\h}_m \oplus \tilde{\h}_n$ written in polar coordinates about a copy of $\H^{m-1} \times \H^{n-1}$.   
In Section \ref{sect:proof of main theorem 2} we prove Theorem \ref{thm:main theorem 2} using the curvature formulas from Theorem \ref{thm: curvature formula polar}.
Finally, in Section \ref{sect:proof of group properties} we prove Theorem \ref{thm: group theoretic properties}.

\bigskip

\noindent \textbf{Acknowledgements:} Cho and Lafont were partially supported by grant number DMS--2407438 and Miller was partially supported by grant number DMS--2405264/2630020 from the National Science Foundation. Minemyer was partially supported by an AMS-Simons Research Enhancement Grant for Primarily Undergraduate Institution Faculty.

In writing this paper, the publicly available version of ChatGPT was used to construct some smoothing functions. For the convenience of the reader, the AI arguments have been isolated to the Appendix. The correctness of the arguments
in the Appendix have been independently verified by the authors. No AI 
was used in the remainder of the paper (Sections 1-5). The final draft was proofread by ChatGPT, and minor corrections were carried out based on the AI report.

\section{Lattices in Isom$\left( \H^m \times \H^n \right)$}\label{sec:background}

\subsection{Generalities on locally symmetric spaces modeled on $\H^m\times \H^n$}
In this subsection, we recall some basic facts about locally symmetric orbifolds $M$ whose universal cover is isometric to $\H^m\times \H^n$, for $m,n\ge 3$.

Throughout we let $q_m$ denote the \emph{standard quadratic form} on $\R^{m+1}$, which is defined by the formula
\begin{equation}\label{eqn:quad}
q_m(x_1,\dots,x_{m+1})=x_1^2+\dots +x_m^2-x_{m+1}^2.
\end{equation}
The hyperboloid model of hyperbolic space is then given by 
\begin{equation}\label{eqn:hyperbolicspace}\H^m=\left\{(x_1,\dots,x_{m+1})\in\R^{m+1}\mid q_m(x_1,\dots,x_{m+1})=-1,~x_{m+1}>0\right\},
\end{equation}
and the orientation preserving group of isometries, $\mathrm{Isom}^+(\H^m)$, of this model is
$$\SO_0(m,1)=\left\{A\in\SL_{m+1}(\R)\mid q_m(Ax)=q_m(x),~\forall x\in\R^{m+1}\right\}^\circ,$$
where the superscript on the righthand side of this equation denotes the identity component.
In particular,
\begin{equation}\label{eqn:fi}
\SO_0(m,1)\times \SO_0(n,1)\cong \mathrm{Isom}^+(\H^m)\times\mathrm{Isom}^+(\H^n)\le_{f.i.}\mathrm{Isom}^+(\H^m\times\H^n),
\end{equation}
where the notation means that this is a subgroup of finite index.
For ease of notation, we let $G_k=\SO_0(k,1)$ throughout the remainder of this section.

Recall that two discrete subgroups $\Gamma,\Gamma'<\mathrm{Isom}^+(\H^m\times\H^n)$ are \emph{commensurable} (in the wide sense) if there exists $g\in \mathrm{Isom}^+(\H^m\times\H^n)$ such that $\Gamma\cap g\Gamma'g^{-1}$ is a finite index subgroup of both $\Gamma$ and $g\Gamma'g^{-1}$. Geometrically this is equivalent to $(\H^m\times\H^n)/\Gamma$ and $(\H^m\times\H^n)/\Gamma'$ having a common finite cover.
The notion of commensurability forms an equivalence relation and we will primarily be interested in commensurability classes of subgroups (equivalently, orbifolds).
As such, Equation \eqref{eqn:fi} implies that we may focus our attention on subgroups $\Gamma<G_m\times G_n$.

We call $\Gamma<G_m\times G_n$ a \emph{lattice} if $\Gamma$ is discrete and the quotient $(G_m\times G_n)/\Gamma$ has finite Haar measure, equivalently, if the orbifold $(\H^m\times\H^n)/\Gamma$ has finite-volume.
If a lattice $\Gamma$ is additionally torsion-free, then the quotient $(\H^m\times\H^n)/\Gamma$ is a finite-volume manifold.
Such lattices are always finitely generated \cite[Cor 13.15]{RagBook} and therefore Selberg's lemma \cite[Lem 8]{Selberg} allows one to find a finite index subgroup $\Gamma'\le \Gamma$ which is torsion-free and consequently a finite-volume manifold $(\H^m\times\H^n)/\Gamma'$ which finitely covers $(\H^m\times\H^n)/\Gamma$. 
Owing to this, for the time being, we postpone the issue of requiring lattices to be torsion-free.

\begin{definition}\label{def:reducible lattices}
We say that a lattice $\Gamma< G_m\times G_n$ is \emph{irreducible} if its projection to each of $G_m$, $G_n$ is dense and \emph{reducible} otherwise.
\end{definition}

\begin{remark}\label{rmk:product lattice finite index}
It is straightforward to construct reducible lattices by letting $\Gamma=\Gamma_1\times \Gamma_2$, where $\Gamma_1<G_m$ and $\Gamma_2<G_n$ are lattices, and indeed all reducible lattices have a finite index subgroup of this form \cite[Ch 5]{RagBook}.
\end{remark}

It is also well-known that irreducible lattices exist in $G_m\times G_n$ if and only if $G_m$ and $G_n$ are isotypic, that is, if the simple factors of their complexifications are isogenous (see for instance \cite[Ch 5F]{DWMBook}). By a straightforward dimension argument, one sees that this is possible if and only if $m=n$ or $m,n\in\{2,3\}$. 
When $m=n$ irreducible lattices are necessarily arithmetic by seminal work of Margulis \cite[Ch IX]{MargulisBook}, as $G_m\times G_n$ is a higher rank Lie group.

\subsection{Arithmetic lattices of simplest type in products of hyperbolic spaces}\label{sec:simpletype}

In this subsection, we describe how to construct families of irreducible arithmetic lattices in $G_m\times G_n$ whenever $m=n\ge 2$. In fact, as it requires minimal extra work, we will describe how to construct infinitely many commensurability classes of arithmetic lattices in $(G_m)^k$ for any $k\in \N$ and therefore finite-volume orbifolds of the form $(\H^m\times\dots\times \H^m)/\Gamma$.
Throughout this subsection, fix a choice of $k\in \N$.

Let $K$ be a totally real number field with $ d=[K:\mathbb{Q}]$ and denote by $\sigma_i:K\to \R$ the Galois embeddings of $K$ for $1\le i\le d$.
We will always require that $d\ge k$ in the sequel.
Let $q$ be a non-degenerate $K$-quadratic form on the $K$-vector space $K^{m+1}$, then we say that the pair $(q,K)$ is \emph{admissible} if the forms $q^{\sigma_1},\dots, q^{\sigma_k}$ each have signature $(m,1)$ over $\mathbb{R}$, and $q^{\sigma_i}$ is definite for any $i>k$.

Given a field $E$ and a quadratic form $q$ on the vector space $E^{m+1}$, let $B_q$ denote the matrix associated to $q$ in the standard basis.  
Then the \emph{special orthogonal group associated to $(E,q)$} is given by
$$\SO_q(E)=\{A\in\SL_{m+1}(E)\mid A^TB_qA=B_q\}.$$
Of particular importance is when $(q,K)$ is an admissible pair, then Sylvester's law of inertia gives that
\begin{equation}\label{eqn:Sylvester}\SO_{q^{\sigma_i}}(\R)\cong \begin{cases}
\SO(m,1),& 1\le i\le k\\
\SO(m+1),& i>k    \end{cases}.
\end{equation}
In particular, we have the embedding
\begin{align*}
\psi:\SO_q(K)&\hookrightarrow \underbrace{\SO(m,1)\times\dots\times \SO(m,1)}_{k\text{ times}}\times \prod_{i>k}\SO(m+1),\\
g&\mapsto \left(\sigma_1(g),\dots,\sigma_k(g),\sigma_{k+1}(g),\dots,\sigma_d(g)\right)
\end{align*}
where we abusively suppress the post-composition with the isomorphisms from Equation \eqref{eqn:Sylvester}.
If $\mathcal{O}_K$ denotes the ring of integers of $K$ and
$$\Lambda=\SO_q(\mathcal{O}_K)=\left\{A\in \SL_{m+1}(\mathcal{O}_K)\mid A^TB_qA=B_q\right\},$$
then $\psi(\Lambda)$ is discrete in $ \SO(m,1)^k\times \SO(m+1)^{d-k}$.
Indeed, this follows from the fact that $\mathcal{O}_K$ is discrete in $\R^d$ via the similar map induced by the Galois embeddings.
An important theorem of Borel--Harish-Chandra \cite{BHC} moreover shows that $\psi(\Lambda)$ is a lattice in this product. 

Letting $\pi$ denote the projection from this product to the $\SO(m,1)^k$ factor, it follows that $\Gamma=\pi(\psi(\Lambda))\cap (G_m)^k$ is a lattice in $(G_m)^k$. 
Indeed, the kernel of $\pi$ is the compact group $\prod_{i>k}\SO(m+1)$. 
It is also straightforward to check from the setup that, when $k\ge 2$, the projection to each $G_m$ factor of $(G_m)^k$ is dense and therefore $\Gamma$ is an irreducible lattice in $(G_m)^k$.

\begin{definition}
We call any lattice in $\mathrm{Isom}^+(\H^m\times\dots\times \H^m)$ that is commensurable with some $\Gamma$ as constructed above an \emph{irreducible arithmetic lattice of simplest type}.
Given such a lattice, we call the pair $(q,K)$ used to construct $\Gamma$ the \emph{arithmetic data} associated to the lattice.
\end{definition}

A couple of facts follow from this definition which we list as remarks.

\begin{remark}\label{rem:commclass}
Given two arithmetic lattices of simplest type $\Gamma$, $\Gamma'$ with associated arithmetic datum $(q,K)$, $(q',K')$, then $\Gamma$ and $\Gamma'$ are commensurable if and only if $K\cong K'$ and $q$ and $q'$ are in the same $K$-similarity class.
The latter means that $q$ is $K$-isometric to $\lambda q'$ for some $\lambda\in K^*$.
In particular, the number field $K$ and the $K$-similarity classes of $q$ is a complete invariant of the commensurability class.
\end{remark}

\begin{remark}\label{rem:anisotropic}
By Godement's compactness criterion \cite[Prop 5.26]{DWMBook} an arithmetic lattice of simplest type $\Gamma$ with associated arithmetic data $(q,K)$ is cocompact if and only if $q$ is anisotropic, that is, if there is no non-zero $v\in K^{m+1}$ for which $q(v)=0$.
Otherwise, we say that $q$ is isotropic.
These notions are commensurability invariants. 
Moreover, the Hasse--Minkowski theorem shows that $q$ is isotropic if and only if $q$ is isotropic over every completion of $K$.
Therefore, whenever $d>k$ the form $q$ is necessarily anisotropic since there is at least one completion associated to a Galois embedding for which $q$ is positive definite (and hence anisotropic).
\end{remark}

For concreteness, we finish the subsection with a few simple families of examples of lattices arising from the above construction.

\begin{example}\label{ex:hyperbolicsimplest}
Let $k=1$ in the discussion above, we now construct two families of arithmetic lattices of simplest type in $\mathrm{Isom}^+(\H^m)$.
Let $K_D=\mathbb{Q}(\sqrt{D})$, where $D\in\N\setminus\{1\}$ is square-free, and define the following two quadratic forms
\begin{align*}
q_m(x_1,\dots,x_{m+1})&=x_1^2+\dots+ x_m^2-x_{m+1}^2,\\
f(x_1,\dots,x_{m+1})&=x_1^2+\dots+ x_m^2-\sqrt{D}x_{m+1}^2,
\end{align*}
where the former is a $\mathbb{Q}$-quadratic form and the latter is a $K_D$-quadratic form.
The pair $(q_m,\mathbb{Q})$ is clearly admissible and if $\sigma_1,\sigma_2$ are the Galois embeddings of $K_D$, one can visibly see that $f^{\sigma_1}$ has signature $(m,1)$ and $f^{\sigma_2}$ has signature $(m+1,0)$, hence the pair $(f,K_D)$ is similarly admissible.
In particular, $$
\Gamma_1=\pi(\psi(\SO_{q_m}(\Z)))\cap G_m,\quad \Gamma_{2,D}=\pi(\psi(\SO_f(\mathcal{O}_{K_D})))\cap G_m,$$ are irreducible arithmetic lattices of simplest type in $G_m$ with arithmetic datum $(q_m,\mathbb{Q})$, $(f,K_D)$ (respectively). 

By the discussion in Remark \ref{rem:anisotropic}, it is straightforward to check that for any $m\ge2$, the quotients $\H^m/\Gamma_1$ are always finite-volume, non-compact orbifolds.
On the other hand, as $d=2>k=1$, Remark \ref{rem:anisotropic} also shows that the orbifolds $\H^m/\Gamma_{2,D}$ are in fact always compact.
Moreover, since $K_D\cong K_{D'}$ if and only if $D=D'$, it follows from Remark \ref{rem:commclass} that for any $m$, as $D$ varies over the set of square-free natural numbers, the latter orbifolds lie in infinitely many distinct commensurability classes.
\end{example}

\begin{example}\label{ex:twohyperbolicsimplest}
We now consider the case $k=2$.
Let $K_D$ as before and $L_{D,D'}=\mathbb{Q}(\sqrt{D},\sqrt{D'})$, where $D, D'\in\N$ are square-free with $D<D'$.
Define the following two quadratic forms
\begin{align*}
h_1(x_1,\dots,x_{m+1})&=(D-\sqrt{D})x_1^2+\dots+ x_m^2-x_{m+1}^2,\\
h_2(x_1,\dots,x_{m+1})&=x_1^2+\dots+ x_m^2-(\sqrt{D}+\sqrt{D'})x_{m+1}^2,
\end{align*}
where the former is a $K_D$-quadratic form and the latter is an $L_{D,D'}$-quadratic form.
Since $h_1$ is $K$-defined, $D-\sqrt{D}$ is not a square in $K_D$, and $D\pm\sqrt{D}>0$, it follows that $h_1^{\sigma_1}$, $h_1^{\sigma_2}$ have signature $(m,1)$, therefore $(h_1,K_D)$ is admissible.
One similarly checks that $(h_2,L_{D,D'})$ is admissible since $\pm\sqrt{D}+\sqrt{D'}>0$ and $\pm\sqrt{D}-\sqrt{D'}<0$ each of which represent the image of $\sqrt{D}+\sqrt{D'}$ under one of the Galois embeddings of the multi-quadratic extension $L_{D,D'}/\mathbb{Q}$.
In particular, $$
\Gamma_{1,D}=\pi(\psi(\SO_{h_1}(\mathcal{O}_{K_D})))\cap (G_m)^2,\quad \Gamma_{2,(D,D')}=\pi(\psi(\SO_{h_2}(\mathcal{O}_{L_{D,D'}})))\cap (G_m)^2,$$ are irreducible arithmetic lattices of simplest type in $G_m\times G_m$ with arithmetic datum $(h_1,K_D)$, $(h_2,L_{D,D'})$ (respectively). 

As before, Remark \ref{rem:commclass} shows that $\Gamma_{1,D}$ is commensurable to $\Gamma_{1,D'}$ if and only if $D=D'$. Similarly, $\Gamma_{2,(D,D')}$, $\Gamma_{2,(E,E')}$ are non-commensurable if the multi-quadratic extensions associated to $(D,D')$, $(E,E')$ differ.
Moreover, $h_1$ is visibly isotropic so by Remark \ref{rem:anisotropic} the finite-volume orbifolds $(\H^m\times \H^m)/\Gamma_{1,D}$ are always non-compact. 
On the other hand, since $[L_{D,D'}:\mathbb{Q}]=4>k=2$ the same remark shows that the finite-volume orbifolds $(\H^m\times \H^m)/\Gamma_{2,(D,D')}$ are always compact.

Letting $D<D'$ vary over the set of all square-free natural numbers we have constructed two infinite families of pairwise non-commensurable, irreducible orbifolds whose universal covers are isometric to $\H^m\times \H^m$; one family which is always closed and the other family which is always finite-volume, non-compact.
To be concrete, if $(D,D')$ varies over $(2,p)$ for $p$ an odd prime then the resulting lattices $\Gamma_{2,(D,D')}$ are pairwise non-commensurable.
\end{example}

\subsection{Totally geodesic subspaces in orbifolds modeled on $\H^m\times \H^n$}\label{sec:geodesic}

Assume that $M$ is an orbifold whose universal cover is isometric to $\H^m\times \H^n$ for $m,n\ge 3$ and similarly $Z$ is an orbifold whose universal cover is isometric to $\H^{m-1}\times \H^{n-1}$.
We say that a map $\iota:Z\to M$ is a \emph{proper totally geodesic immersion} if $\iota$ is proper and any of the lifts $\widetilde{\iota}:\H^{m-1}\times\H^{n-1}\to \H^{m}\times\H^{n}$ is a totally geodesic embedding.
In accordance, we say that $S=\iota(Z)\subset M$ is an \emph{immersed totally geodesic subspace}. When $\iota$ is an embedding, we call $S$ \emph{embedded}, and when $Z$ is a manifold, we indicate it as such.

If $\Gamma=\pi_1(M)$, properness is equivalent to requiring that the subgroup $$\Pi=\mathrm{stab}_{\Gamma}\left(\widetilde{\iota}\left(\H^{m-1}\times\H^{n-1}\right)\right)<\Gamma,$$
which stabilizes the lift $\widetilde{\iota}$, is a lattice in $\mathrm{stab}_{\mathrm{Isom}^+(\H^m\times \H^n)}\left(\widetilde{\iota}\left(\H^{m-1}\times\H^{n-1}\right)\right)$. 
Let
\begin{align}
j:G_{m-1}\times G_{n-1}&\hookrightarrow G_m\times G_n.\label{eqn:block}\\
(A,B)&\mapsto \left(\begin{pmatrix}
1&0\\
0&A\end{pmatrix},\begin{pmatrix}
1&0\\
0&B\end{pmatrix}\right)\nonumber
\end{align}
be the block diagonal embedding of $G_{m-1}\times G_{n-1}$ into $G_m\times G_n$.
Then $j(G_{m-1}\times G_{n-1})$ stabilizes $\H^{m-1}_{\mathrm{std}}\times\H^{n-1}_{\mathrm{std}}$, where $\H^{m-1}_{\mathrm{std}}$ denotes the codimension $1$ isometrically embedded copy of $\H^{m-1}$ in $\H^m$ given by taking $x_1=0$ in Equation \eqref{eqn:hyperbolicspace} (similarly for $\H^n$).
Moreover, $j(G_{m-1}\times G_{n-1})$ is finite index in the full stabilizer of $\H^{m-1}_{\mathrm{std}}\times\H^{n-1}_{\mathrm{std}}$ in $\mathrm{Isom}^+(\H^m\times \H^n)$.

Since $G_m$ (resp. $G_n$) acts transitively on codimension $1$ hyperbolic subspaces of $\H^m$ (resp. $\H^n$), it follows that there exists $g\in G_m\times G_n$ for which the following chain of finite index inclusions hold
\begin{align*}
j(G_{m-1}\times G_{n-1})^{g}&\le _{f.i} \mathrm{stab}_{G_m\times G_n}\left(\widetilde{\iota}\left(\H^{m-1}\times\H^{n-1}\right)\right)\\
&=\left(\left(\mathrm{O}(1)\times \mathrm{O}(m-1,1)\right)\times\left(\mathrm{O}(1)\times\mathrm{O}(n-1,1)\right)\right)^g\cap \left(G_m\times G_n\right)\\
&\le_{f.i}\mathrm{stab}_{\mathrm{Isom}^+(\H^m\times \H^n)}\left(\widetilde{\iota}\left(\H^{m-1}\times\H^{n-1}\right)\right),\end{align*}
where the superscript denotes conjugation.
Since the property of being a lattice is invariant under passage to finite index, we therefore conclude that constructing immersed totally geodesic submanifolds modeled on $\H^{m-1}\times \H^{n-1}$ is equivalent to finding elements $g\in G_m\times G_n$ for which 
$$\Delta=j(G_{m-1}\times G_{n-1})^g\cap \Gamma<G_m\times G_n,$$
is a lattice in $j(G_{m-1}\times G_{n-1})^g$.

\subsection{Examples of manifolds satisfying the hypothesis of Theorems \ref{thm:main theorem 1} and \ref{thm:main theorem 2}}

Using the previous subsections, for each pair $m,n\ge 3$, we now provide infinitely many commensurability classes of manifolds to which Theorems \ref{thm:main theorem 2} and \ref{thm:main theorem 1} apply.

\begin{theorem}\label{thm:construction}
For every $m,n\in\mathbb{N}$ such that $m,n\ge 3$, there exists infinitely many commensurability classes of finite-volume manifolds with universal cover isometric to $\H^m\times \H^n$ and containing an embedded compact totally geodesic submanifold with universal cover isometric to $\H^{m-1}\times \H^{n-1}$. 
Moreover, when $m=n$, these manifolds may be taken to be either reducible or irreducible.
\end{theorem}

We do not claim novelty in the statement or the proof of Theorem \ref{thm:construction}; indeed, this is well-known to the reader fluent in arithmetic constructions.
However we attempt, as we have above, to give the relevant details for readers who have not seen such constructions in the literature before.

\begin{proof}[Proof of Theorem \ref{thm:construction} in the reducible case]
Let $\Lambda_1<G_m$ be a fixed cocompact arithmetic lattice of simplest type and $\Lambda_{2,i}<G_n$ be an infinite family of pairwise non-commensurable cocompact arithmetic lattices of simplest type. For example, one may concretely take $\Lambda_1=\Gamma_{f,2}$, $\Lambda_{2,i}=\Gamma_{f,D_i}$ as $D_i$ varies over the set of square-free natural numbers from Example \ref{ex:hyperbolicsimplest}.
It is well-known that orbifolds arising as quotients by arithmetic lattices of simplest type in $\mathrm{Isom}^+(\H^m)$ always have immersed codimension $1$ totally geodesic subspaces (in fact infinitely many) and, moreover, that given any fixed compact immersed codimension $1$ totally geodesic subspace $S_1\subset \H^m/\Lambda_1$ (resp. $S_{2,i}\subset \H^n/\Lambda_{2,i}$) there is a finite index subgroup $\Gamma_1\le_{f.i.}\Lambda_1$ (resp. $\Gamma_{2,i}\le_{f.i.}\Lambda_{2,i}$) for which the following hold:
\begin{enumerate}
    \item $\Gamma_1$ (resp. $\Gamma_{2,i}$) is torsion-free, and
    \item there is a lift $Y_1$ of $S_1$ to $\H^m/\Gamma_1$ (resp. $Y_{2,i}$ of $S_{2,i}$ to $\H^n/\Gamma_{2,i}$) which is embedded.
    \end{enumerate}
Indeed, the former follows from Selberg's lemma and a proof of the latter can be found in \cite[Lem 1.7]{BHW} (alternatively, we recount the argument in the proof of the irreducible case below).
Consequently, $\Gamma_i=\Gamma_1\times \Gamma_{2,i}$ is a reducible lattice in $G_m\times G_n$ with the property that $S=Y_1\times Y_{2,i}$ is a closed embedded totally geodesic submanifold of the closed manifold $M_i=(\H^m\times \H^n)/\Gamma_i$.
Moreover, that the $M_i$ are pairwise non-commensurable follows from the pairwise non-commensurability of $\Gamma_{2,i}$.
\end{proof}

\begin{proof}[Proof of Theorem \ref{thm:construction} in the irreducible case.]
Let $(q,K)$ be admissible for $k=2$ with $d=[K:\mathbb{Q}]> 2$.
Up to $K$-similarity, we can and do assume that $q$ is diagonal, that is, of the following form
$$q(x_1,\dots,x_{m+1})=a_1x_1^2+a_2x_2^2+\dots +a_{m+1}x_{m+1}^2,$$
where $a_i\in K^*$ and moreover, we may further assume that $\sigma_i(a_1)>0$ for all $1\le i\le d=[K:\mathbb{Q}]$ due to the signature condition.
We let $$\Lambda=\pi(\psi(\SO_q(\mathcal{O}_K)))\cap\left(G_m\times G_m\right) <G_m\times G_m,$$
be the corresponding lattice, as described in Section \ref{sec:simpletype}.

Let $W\subset K^{m+1}$ be the $K$-hyperplane defined by setting $x_1=0$ and note that, by the condition on $a_1$, the restriction $q_0=q\vert_W$ is a non-degenerate quadratic form on $W$. Moreover, if 
$$\Xi=\pi(\psi(\SO_{q_0}(\mathcal{O}_K)))\cap\left(G_{m-1}\times G_{m-1}\right),$$
then it follows that $j(\Xi)<j(G_{m-1}\times G_{m-1})<G_m\times G_m$ and moreover, by construction, that $j(\Xi)$ is also an irreducible lattice in $j(G_{m-1}\times G_{m-1})$.
Therefore following the discussion of Section \ref{sec:geodesic} there is a subgroup $\Pi<\Lambda$ for which $\left(\H^{m-1}_{\mathrm{std}}\times \H^{m-1}_{\mathrm{std}}\right)/\Pi$ is an immersed totally geodesic subspace of $\left(\H^m\times \H^m\right)/\Lambda$.
By Selberg's lemma, we may always find a torsion-free finite index subgroup $\Upsilon\le_{f.i.}\Lambda$ and therefore $\left(\H^m\times \H^m\right)/\Upsilon$ is a manifold with the requisite immersed totally geodesic submanifold $\left(\H^{m-1}_{\mathrm{std}}\times \H^{m-1}_{\mathrm{std}}\right)/(\Pi\cap\Upsilon)$.

Finally, we upgrade this immersion to an embedding as follows.
The subspace $\left(\H^{m-1}_{\mathrm{std}}\times \H^{m-1}_{\mathrm{std}}\right)/\left(\Pi\cap\Upsilon\right)$ is not embedded if and only if the set
$$\mathcal{B}=\left\{\upsilon\in\Upsilon\setminus \mathrm{stab}_{\Upsilon}\left(\H^{m-1}_{\mathrm{std}}\times \H^{m-1}_{\mathrm{std}}\right)~\middle|~\upsilon\left(\H^{m-1}_{\mathrm{std}}\times \H^{m-1}_{\mathrm{std}}\right)\cap \left(\H^{m-1}_{\mathrm{std}}\times \H^{m-1}_{\mathrm{std}}\right)\neq\emptyset\right\},$$
is non-empty.
Letting 
$V=\left(\mathrm{O}(1)\times\mathrm{O}(m-1,1)\right)^2$, which is a real algebraic subset of $\mathrm{GL}_{2(m+1)}(\R)$ for which $\mathrm{stab}_{G_m\times G_m}\left(\H^{m-1}_{\mathrm{std}}\times \H^{m-1}_{\mathrm{std}}\right)\subset V$, we have the following chain of 
inclusions
$$\Upsilon <G_m\times G_m<\mathrm{O}(m,1)\times\mathrm{O}(m,1)<\mathrm{GL}_{2(m+1)}(\R).$$
As $\Upsilon$ is finitely generated, the subgroup
$$\mathrm{stab}_{\Upsilon}(\H^{m-1}_{\mathrm{std}}\times\H^{m-1}_{\mathrm{std}})=V\cap\Upsilon=\Pi\cap\Upsilon,$$
is separable by \cite[Lem 1.8]{BHW}.
Since $q_0$ is definite at a real place, the quotient $(\H^{m-1}_{\mathrm{std}}\times\H^{m-1}_{\mathrm{std}})/(\Pi\cap \Upsilon)$ is compact.
The standard compact-lifting criterion for separable subgroups (see for instance \cite[Lem 1.7]{BHW}) then gives a finite-index subgroup $\Gamma\le_{f.i.}\Upsilon$ containing $\Pi\cap\Upsilon$ for which $\mathcal{B}\cap\Gamma=\emptyset$.
Therefore $M=\left(\H^{m}\times \H^{m}\right)/\Gamma$ is a manifold containing an embedded totally geodesic submanifold of the form $S=\left(\H^{m-1}\times \H^{m-1}\right)/(\Gamma\cap \Pi)$.
An infinite family of commensurability classes of lattices satisfying the criterion at the beginning of the proof is given in Example \ref{ex:twohyperbolicsimplest}.
\end{proof}

\begin{remark}
    For a lattice in the product $G_m\times G_m$, there is an alternate geometric criterion for reducibility: one can take the orbit of an arbitrary point $p\in \mathbb H^m \times \mathbb H^m$ under the action of the lattice, and look at the projection of this orbit onto one of the $\mathbb H^m$ factors. Then the lattice is reducible if and only if the image under the projection map is a discrete subset of $\mathbb{H}^m$. This immediately implies that, in the case where the ambient arithmetic lattice $\Gamma$ is reducible, the lattice $\Lambda$ associated to the submanifold will also have to be reducible. 
\end{remark}

\section{Curvature formulas}\label{Section:Prelim}

In this section we first provide details for a construction that was quickly outlined in the second and third paragraphs of Section 3 of \cite{fornari1990ramified}.  
We provide full details here for two reasons.
The first is because this construction is critical to our proof for Theorem \ref{thm:main theorem 1} (1) (see Section \ref{sect:product case}) and serves to illustrate why this argument only works for product lattices. The second reason is that this computation serves as an introduction to more complicated techniques used later in the paper.
We then develop curvature formulas for a more sophisticated warped product metric (Theorem \ref{thm: curvature formula polar}) that will be used in Section \ref{sect:proof of main theorem 2}.

\subsection{Warped product metrics on $\H^n$ about $\H^{n-1}$} Consider the manifold $\R^n \cong \R^{n-1} \times (- \infty, \infty)$.  
If we endow $\R^{n-1}$ with the hyperbolic metric $\tilde{\h}_{n-1}$ and let $x_n$ denote the variable in the last coordinate, then it is well-known (see for example \cite{belegradek2007rigidity}) that the hyperbolic metric $\tilde{\h}_n$ on $\R^n$ can be written as
    \begin{equation}\label{eqn:projection metric}
    \tilde{\h}_n= \cosh^2(x_n)\tilde{\h}_{n-1}+dx_{n}^2.
    \end{equation}

Now, let $\sigma_n : \R \to (0, \infty)$ be a positive, real-valued function of $x_n$.  
Replacing $\cosh(x_n)$ with $\sigma_n(x_n)$ in Equation \eqref{eqn:projection metric} yields the {\it warped-product metric}
    \begin{equation}\label{eqn:1d warped product metric}
    \tilde{g}_n = \sigma_n^2(x_n) \tilde{\h}_{n-1} + dx_n^2.
    \end{equation}
The components of the sectional curvature tensor of $\tilde{g}_n$ will, of course, be functions of $\sigma_n$.  
Formulas for the sectional curvature tensor of $\tilde{g}_n$ are again well-known (see for example \cite[Eqns 3.1, 3.2]{belegradek2007rigidity}).  
We list the relevant formulas after establishing an appropriate frame.  
All calculations in the following paragraph are with respect to the hyperbolic metric $\tilde{\h}_n$.  

Let $\varphi: \R^n \to \R^{n-1} \times \{ 0 \}$ be the orthogonal projection.  
Let $q \in \R^n$ and $p = \varphi(q)$.  
Choose an orthonormal frame $(
\check{X}_1, \hdots , \check{X}_{n-1})$ about $p$ for $T_p \R^{n-1}$.
Since $\varphi: \R^n \to \R^{n-1}$ is a submersion, $d \varphi : T_q \R^n \to T_p \R^{n-1}$ is an epimorphism.  
Letting $\mathcal{H}$ denote the orthogonal complement to $Y_n := \partial / \partial x_n$ in $T_q \R^n$, we have that $d \varphi\vert_\mathcal{H} : \mathcal{H} \to T_p \R^{n-1}$ is an isomorphism.  
Letting $X_i = d \varphi\vert_\mathcal{H}^{-1} ( \check{X}_i ) \in \mathcal{H}$ defines an orthogonal frame $(X_1, \hdots, X_{n-1})$ for $\mathcal{H}$ about $q$.  
This frame is {\it not} orthonormal since $\tilde{\h}_n(X_i, X_i) = \cosh^2(x_n)$ for each $i$.  
So the collection of vectors
    \begin{equation}
    Y_i = \frac{1}{\cosh(x_n)} X_i \quad \text{for } 1 \leq i \leq n-1, \qquad Y_n = \frac{\partial}{\partial x_n}
    \end{equation}
forms an orthonormal frame for $T_q \R^n$ about $q$ with respect to the metric $\tilde{\h}_n$.  

The corresponding orthonormal frame for $\tilde{g}_n$ is
    \begin{equation}\label{eqn:orthonormal frame}
    Y_i = \frac{1}{\sigma_n} X_i \quad \text{for } 1 \leq i \leq n-1, \qquad Y_n = \frac{\partial}{\partial x_n}.
    \end{equation}

Let $K_n(A, B)$ denote the sectional curvature of the $2$-plane spanned by the vectors $A$ and $B$ with respect to $\tilde{g}_n$.  
Formulas for the sectional curvatures of $\tilde{g}_n$ with respect to the frame in Equation \eqref{eqn:orthonormal frame} are
    \begin{align}
        &K_n(Y_i, Y_j) = - \frac{1}{\sigma_n^2} - \left( \frac{\sigma_n'}{\sigma_n} \right)^2  \quad \text{ for } 1 \leq i, j \leq n-1 \text{ and }  i \neq j. \label{eqn: curvature formula 1} \\
        &K_n (Y_i, Y_n) = - \frac{\sigma_n''}{\sigma_n} \quad \text{ for } 1 \leq i \leq n-1.  \label{eqn: curvature formula 2}
    \end{align}
Additionally, there are no mixed terms.  
The frame in Equation \eqref{eqn:orthonormal frame} diagonalizes the curvature operator corresponding to $\tilde{g}_n$.

\begin{proposition}\label{prop:interpolating metric 1}
Let $\eta > 0$.  Then there exist $\e_1$ and $\e_2$ satisfying $0 < \e_1, \e_2 < (1/2) \eta$ and a non-positively curved metric $g$ on $\R^{n-1} \times [- \eta, \eta]$ which satisfies 
    \begin{align*}
    &\tilde{g} = \tilde{\h}_{n-1} + dx_n^2 \quad \text{ over } \, \R^{n-1} \times (-\e_1, \e_1), \\
    &\tilde{g} = \tilde{\h}_n  \hskip 48pt \text{ over } \, \R^{n-1} \times \left( [- \eta , - \eta + \e_2) \cup (\eta - \e_2, \eta] \right).
    \end{align*}
\end{proposition}

Proposition \ref{prop:interpolating metric 1} is easily seen to be true for any $\e_1, \e_2$ which satisfy the given inequality, however, in what follows we only need the weaker statement above. 

\begin{proof}
Let $\e_1 = \e_2 = (1/8) \eta$.  
Define the continuous, piecewise-smooth function $\hat{\sigma}_n$ over $[0, \eta]$ by
    \[ \hat{\sigma}_n =  \begin{cases} 
      	 1, & \text{if } 0 \leq x \leq \frac{1}{4} \eta \\
     	 1 + \frac{\cosh \left( \frac{3}{4} \eta \right) - 1}{\frac{1}{2} \eta} \left( x - \frac{1}{4} \eta \right), & \text{if } \frac{1}{4} \eta \leq x \leq \frac{3}{4} \eta \\
      	\cosh(x), & \text{if } \frac{3}{4} \eta \leq x \leq \eta
   	\end{cases}
	\]
The function $\hat{\sigma}_n$ has singularities at $(1/4) \eta$ and $(3/4) \eta$, but is easily seen to satisfy the hypothesis of Proposition \ref{lem:smoothing lemma} in the Appendix.  
This statement allows us to replace $\hat{\sigma}_n$ by a smooth function $\sigma_n$ extended in a symmetric manner to the interval $[-\eta , \eta]$.  
Then the associated metric $\tilde g_n$ given by the formula in equation (\ref{eqn:1d warped product metric}) has the desired properties, per our earlier calculations for the curvature (see equations (\ref{eqn: curvature formula 1}) and (\ref{eqn: curvature formula 2})) and the properties of the smoothing from Proposition \ref{lem:smoothing lemma}.   
\end{proof}

\begin{remark}\label{rmk:tubular-nbhd-flattening}
    The construction in Proposition \ref{prop:interpolating metric 1} is invariant under isometries of $\h_{n-1}$ acting on the $\R^{n-1}$-factor, as well as under the involution $x_n \mapsto -x_n$ on the last coordinate. It follows that if $S$ is a closed embedded totally geodesic codimension one submanifold inside a hyperbolic manifold $M$, then one can modify the metric in the $\eta$-neighborhood of $S$ to be locally isometric to $\H^{n-1}\times \R$ near $S$.
\end{remark}

\begin{remark}\label{rmk:bi-Lipschitz}
    We can use the identity map on $\R^{n-1} \times [-\eta, \eta]$ to compare the $\tilde g$ metric to the hyperbolic metric $\tilde \h_n$ in Proposition \ref{prop:interpolating metric 1}. Note that we have the uniform bound $\frac{1}{\cosh(\eta)}\leq \frac{\tilde g(v)}{\tilde \h_n(v)}\leq \cosh(\eta)$ for all tangent vectors, so the two metrics are $C$-bilipschitz to each other, where the bilipschitz constant $C:= \cosh(\eta)$ tends to $1$ as $\eta$ tends to zero. This implies pointwise bilipschitz control of the respective volume forms. In particular, for any $\epsilon >0$, we can choose $\eta$ close enough to zero to ensure that the perturbed metric discussed in Remark \ref{rmk:tubular-nbhd-flattening} has volume $\epsilon$-close to the volume of the hyperbolic metric. 
\end{remark}

\subsection{Local flattening in product lattice case}  

Let us now discuss an application to the case where $\Gamma < \text{Isom}(\H^m \times \H^n)$ is a product lattice, and $S$ is connected.
So there exist $\Gamma_1 < \text{Isom}(\H^m)$ and $\Gamma_2 < \text{Isom}(\H^n)$ such that $\Gamma = \Gamma_1 \times \Gamma_2$. Then $M$ splits isometrically as a direct product $M=X_1 \times X_2$, where $X_1 = \H^m / \Gamma_1$ and $X_2 = \H^n / \Gamma_2$.  
From our hypothesis, the universal cover of our codimension two submanifold $S$ is isometric to $\H^{m-1} \times \H^{n-1}$. This implies $S$ naturally splits as a product $Y_1 \times Y_2$, where $Y_i \subset X_i$ are a pair of connected, embedded, codimension one totally geodesic submanifolds.
In what follows, we use $e_2$ to denote the Euclidean metric on $\R^2$.

\begin{proposition}\label{prop:product metrics}
Let $\Gamma < \Isom^+(\H^m \times \H^n)$ be a product lattice, with quotient manifold $M = (\H^m \times \H^n) / \Gamma$ equipped with the locally symmetric metric $\h_m \oplus \h_n$. Let $S \subset M$ be a connected, embedded, totally geodesic codimension two submanifold, with the property that the universal cover $\tilde S$ is isometric to 
$\H^{m-1} \times \H^{n-1}$. Then for any $\epsilon >0$ there exists a metric $g$ on $M$ having the following properties:
\begin{itemize}
    \item $g$ is non-positively curved,
     \item $g$ is isometric to $\h_{m-1} \oplus \h_{n-1} \oplus e_2$ inside a sufficiently small neighborhood of $S$,
    \item $|\Vol(M, g) -\Vol(M, \h_m \oplus \h_n)|< \epsilon$.
\end{itemize}

\end{proposition}

\begin{proof}
    Consider the codimension one submanifold $Y_1 \subset X_1$, which is a connected, embedded, totally geodesic submanifold inside $X_1$. Applying Proposition \ref{prop:interpolating metric 1} and Remark \ref{rmk:tubular-nbhd-flattening}, $X_1$ supports a metric $g_1$ which is non-positively curved, and locally a flat product in a small neighborhood of $Y_1$. Moreover the metric can be chosen to have volume as close to the hyperbolic one as desired. Similarly, the other factor $X_2$ has a metric $g_2$ which is locally a product near $Y_2$. The metric $g$ is then defined to be the product Riemannian metric $g:=g_1 \oplus g_2$ on $M=X_1\times X_2$. It is clear from Proposition \ref{prop:interpolating metric 1} and Remark \ref{rmk:tubular-nbhd-flattening} that $g$ has the desired properties listed in the Proposition.
\end{proof}

\begin{remark}\label{rmk:not-connected-S}
    Note that $S$ connected plays a minor role in our proof. If $S$ has multiple connected components, the main concern is the possibility that the codimension one submanifolds obtained by projecting onto one of the factors $X_i$ could have intersections. As long as the projections of the individual components of $S$ are pairwise disjoint, the same argument can be carried out.
\end{remark}

\begin{remark}\label{rmk:only for product lattices}
    Observe that our proof crucially relies on the global product structure on the manifold $M= X_1 \times X_2$. The metric ${g}_i$ differs from ${\h}_i$ on some strip about $Y_i$. As a result, the metrics $g$ and $\h_m \oplus \h_n$ differ in strips about the subset $(Y_1 \times X_2) \cup (X_1 \times Y_2)$ -- a much larger set than the submanifold $S= Y_1\times Y_2$. In other words, the metric $g$ is {\bf not} a local perturbation of the metric $\h_m \oplus \h_n$ near the totally geodesic codimension two submanifold $S$, but rather a global perturbation exploiting the product structure. 

    If we consider instead irreducible lattices, then we no longer have a global product structure to exploit. Instead, we have to consider perturbations of the metric concentrated near the submanifold $S$. This requires the more sophisticated warped product constructions discussed in the next section.

    One can try to extend this result to the case of {\it reducible} lattices, by arguing on a finite cover which is a product. But we then run into the issue that the lift of $S$ might not be connected (see the discussion above in Remark \ref{rmk:not-connected-S}).
\end{remark}

\subsection{Polar coordinates, the warped-product metric, and curvature formulas} 
There are several easy arguments to prove Theorem \ref{thm:main theorem 1} $(1)$ using Proposition \ref{prop:product metrics}.
However, to prove Theorem \ref{thm:main theorem 2} we will need the more sophisticated curvature formulas in Theorem \ref{thm: curvature formula polar} below.  

Consider the product $\R^m \times \R^n \cong \R^{m-1} \times (- \infty, \infty) \times \R^{n-1} \times (- \infty , \infty)$.  
By Equation \eqref{eqn:projection metric}, we have that the product hyperbolic metric can be written as
    \begin{equation*}
    \tilde{\h}_m \oplus \tilde{\h}_n = \cosh^2(x_m) \tilde{\h}_{m-1} + dx_m^2 + \cosh^2(y_n)\tilde{\h}_{n-1} + dy_n^2.
    \end{equation*}
If we switch to polar coordinates  
    \begin{equation*}
    x_m = r \cos (\theta), \quad y_n = r \sin(\theta),
    \end{equation*}
for $r \in [0, \infty)$ and $\theta \in \S^1$, then the product hyperbolic metric becomes
    \begin{equation}\label{eqn:hyperbolic polar coordinates metric}
    \tilde{\h}_m \oplus \tilde{\h}_n = \cosh^2(r \cos(\theta)) \tilde{\h}_{m-1} + \cosh^2(r \sin(\theta)) \tilde{\h}_{n-1} + r^2 d \theta^2 + dr^2.
    \end{equation}
The associated warped-product metric is
\begin{equation}\label{eqn:polar coordinate warped product metric}
    \tilde{\l} = \a^2(r, \th) \tilde{\h}_{m-1} + \b^2(r, \th) \tilde{\h}_{n-1} + f^2(r) d \th^2 + dr^2,
    \end{equation}
where $\a$ and $\b$ are positive functions of $r$ and $\th$ and $f$ is a positive function of $r$.  
In Theorem \ref{thm: curvature formula polar} we work out formulas for the components of the sectional curvature tensor of $\l$.  
These formulas will, of course, be functions of $\a$, $\b$, and $f$.   

As in Section \ref{Section:Prelim}, we must first define an orthonormal frame for $\lambda$. 
For the hyperbolic metrics in $\l$, we use the half-space model and write 
\begin{align}
&\h_{m-1}=\frac{1}{x^2_{m-1}}dx_1^2+\dots+\frac{1}{x^2_{m-1}}dx_{m-1}^2,\\
&\h_{n-1}=\frac{1}{y^2_{n-1}}dy_1^2+\dots+\frac{1}{y^2_{n-1}}dy_{n-1}^2.
\end{align}

Let $X'_i=\frac{\partial}{\partial x_i}$, for $1\leq i\leq m-1$ and $Y'_j=\frac{\partial}{\partial y_j}$, for $1\leq j\leq n-1$. Then $(X'_1,\dots,X'_{m-1})$ and $(Y'_1,\dots,Y'_{n-1})$ will be an orthogonal frame for the orthogonal complement to $\big(\frac{\partial}{\partial \theta} ,\frac{\partial}{\partial r}\big)$ defined in an identical manner as in Subsection \ref{Section:Prelim}.1. 

We define the orthonormal frame $(X_1,\dots,X_{m
-1},Y_1,\dots,Y_{n-1}, Z_\theta,Z_r)$ by 
\begin{align}
    &X_i=\frac{x_{m-1}}{\a(r,\theta)}X'_i \quad \text{for } 1 \leq i \leq m-1, \label{eqn:ON basis 1} \\
    &Y_i=\frac{y_{n-1}}{\b(r,\theta)}Y'_{i} \quad \text{for } 1 \leq i \leq n-1,\\
    &Z_{\theta}=\frac{1}{f(r)}\frac{\partial}{\partial \th},\\
    &Z_{r}=\frac{\partial}{\partial r}. \label{eqn:ON basis 4}
\end{align}

In what follows, for the functions $\a(r, \th)$ and $\b(r, \th)$, we will use the standard notation
    $$ \a_r = \frac{\partial \a}{\partial r}, \quad \a_\th = \frac{\partial \a}{\partial \th}, $$
and similarly for $\beta$ and for second order partial derivatives.
It is important to note that subscripts only stand for partial derivatives when they are used with the functions $\a$ and $\b$.

\begin{theorem}\label{thm: curvature formula polar}
    Let $\tilde{\l} = \a^2(r, \th) \tilde{\h}_{m-1} + \b^2(r, \th) \tilde{\h}_{n-1} + f^2(r) d \th^2 + dr^2$ be the metric on $\R^{m-1} \times \R^{n-1} \times \S^1 \times [0, \infty)$ defined in Equation \eqref{eqn:polar coordinate warped product metric}. 
    Up to the symmetries of the curvature tensor, the non-zero components of the curvature tensor are:
    \begin{align}
        &K_\l(X_i, X_j) =  -\frac{1}{\a^2}\left(1+\a_r^2+\frac{\a_\th^2}{ f^2}\right) \quad\text{where} \ i\neq j, \hskip 61pt 
        \label{eqn:curvature component}   \\
        &K_\l(Y_i,Y_j) = -\frac{1}{\b^2} \left( 1 + \b_r^2 + \frac{\b_\th^2}{f^2} \right) \quad\text{where} \ i\neq j, \hskip 67pt  \label{eqn:curvature component 3}
        \\ 
        &K_\l(X_i, Z_\theta) = -\frac{1}{\a f} \left( \a_r f' +\frac{\a_{\th\th}}{f} \right), \hskip 104 pt  \label{eqn:curvature component 4} \\   
        &K_\l(Y_i, Z_\theta) = -\frac{1}{\b f} \left( \b_r f' + \frac{\b_{\th\th}}{f} \right), \hskip 108 pt  \label{eqn:curvature component 5} \\
        &K_\l(X_i, Z_r) = -\frac{\a_{rr}}{\a}, \hskip 168 pt  \label{eqn:curvature component 6}  \\
        &K_\l(Y_i, Z_r) = -\frac{\b_{rr}}{\b}, \hskip 172 pt  \label{eqn:curvature component 7}  \\
        &K_\l(X_i, Y_j) = -\frac{1}{\a \b} \left( \a_r \b_r + \frac{\a_\th\b_\th}{f^2} \right), \hskip 119pt  \label{eqn:curvature component 8}  \\
        &K_\l(Z_r, Z_\theta) = -\frac{f''}{f}. \hskip 196pt  \label{eqn:curvature component 2} 
    \end{align}    
    For non-zero mixed terms, we only have
    \begin{align}
    &<R_\l(X_i,Z_r)X_i,Z_\th>_\l=
           \frac{1}{\a f} \left( \frac{\a_\th f'}{f}- \a_{r\th} \right), \hskip 40pt  \label{eqn:mixed term} \\
    &<R_\l(Y_i,Z_r)Y_i,Z_\th>_\l=
           \frac{1}{\b f} \left( \frac{\b_\th f'}{f}-\b_{r\th} \right). \hskip 48pt  \notag
    \end{align}
\end{theorem}

\begin{remark}
If one inserts the functions $\a = \cosh(r \cos \th)$, $\b = \cosh(r \sin \th)$, and $f = r$ into the metric in Equation \eqref{eqn:polar coordinate warped product metric}, then the metric $\tilde{\l}$ is the product hyperbolic metric $\tilde{\h}_m + \tilde{\h}_n$.  
The curvatures for the metric $\tilde{\h}_m + \tilde{\h}_n$ are given by:  

\begin{align*}
    K_\l(X_i, X_j) = -1 \hskip 124pt & K_\l(Y_i, Y_j) = -1 \\
     K_\l(X_i, Z_\th) = -\sin^2(\th) \hskip 100pt & K_\l(Y_i, Z_\th) = -\cos^2(\th) \\
    K_\l(X_i, Z_r) = -\cos^2(\th) \hskip 100pt & K_\l(Y_i, Z_r) = -\sin^2(\th) \\
     K_\l(X_i, Y_j) = 0 \hskip 136pt & K_\l(Z_r, Z_\th) = 0 \\
    <R_\l(X_i,Z_r)X_i,Z_\th>_\l = \sin(\th) \cos(\th) \\
    <R_\l(Y_i, Z_r) Y_i,Z_\th>_\l= -\sin(\th) \cos(\th) \\
\end{align*}
\noindent One can check that the curvature formulas listed above in Theorem \ref{thm: curvature formula polar}, evaluated with the specific values of $\a = \cosh(r \cos \th)$, $\b = \cosh(r \sin \th)$, and $f = r$, reduce to the expressions above. 
\end{remark}

\begin{proof}

The authors calculated the above formulas using two different methods and obtained the same results each way.  
We will describe both methods, and illustrate their use by giving details for a few of the calculations. The first method  will be used to calculate one of the formulas in Equation \eqref{eqn:curvature component}, and the second method to compute one of the formulas in Equation \eqref{eqn:mixed term}.

The first method is to calculate everything directly using local coordinates and Christoffel symbols.
Recall that, for a metric $g$ and its coefficient matrix $[g_{ij}]$, the formulas for the curvature tensor with orthonormal frame $(X_1,\dots,X_n)$ are 
\begin{equation}
    \la R(X_i,X_j)X_k,X_s \ra_\l =R_{ijks}=\sum_l R^l_{ijk}g_{ls},\label{eqn:tensor1}
\end{equation}
where
\begin{equation}
    R^s_{ijk}=\sum_l \Gamma^l_{ik}\Gamma^s_{jl}-\sum_l \Gamma^l_{jk}\Gamma^s_{il}+\frac{\partial}{\partial x_j}\Gamma^s_{ik}-\frac{\partial}{\partial x_i}\Gamma^s_{jk},\label{eqn:tensor2}  
\end{equation}
and
\begin{equation}
    \Gamma^m_{ij}=\frac{1}{2}\sum_l \left( \frac{\partial}{\partial x_{i}}g_{jl}+\frac{\partial}{\partial x_{j}}g_{li}-\frac{\partial}{\partial x_{l}}g_{ij} \right) g^{lm}.\label{eqn:christoffel}
\end{equation}
Note that, in Equation \eqref{eqn:christoffel}, we are using the standard notation $[g^{ij}]=[g_{ij}]^{-1}$ for the entries in the inverse of the matrix representing the metric.

Now, for the metric $\tilde{\l} = \a^2(r, \th) \tilde{\h}_{m-1} + \b^2(r, \th) \tilde{\h}_{n-1} + f^2(r) d \th^2 + dr^2$ and its orthogonal frame $(X'_1,\dots,X'_{m-1},Y'_1,\dots,Y'_{n-1},\frac{\partial}{\partial\th},\frac{\partial}{\partial r})$, the $(m+n)\times(m+n)$-matrix of the coefficient functions for the metric $\l$ is the diagonal matrix
    \begin{align}\label{eqn:matrix}           
    \left[ \tilde{\l}_{ij} \right]= \begin{bmatrix}
              \frac{\a^2}{x^2_{m-1}} & 0 & \hdots & &  & & \hdots & 0 \\
               0 & \ddots &  &  &  &  & & \vdots \\
               \vdots &  & \frac{\a^2}{x^2_{m-1}} &  &  &  & &\\
               &  &  & \frac{\b^2}{y^2_{n-1}} &  &  &&\\
               & & & & \ddots & & &\\
               &  &  &  & & \frac{\b^2}{y^2_{n-1}} && \\
               \vdots &  & & & &&f^2 & 0 \\
               0 & \hdots  & & & &&0 & 1
              \end{bmatrix}.
    \end{align}

In this computation, we label $m+n$ indices not just from $1$ to $m+n$, but label $x_1,\dots,x_{m-1}$ for the first $m-1$ indices, $y_1,\dots,y_{n-1}$ for the next $n-1$ indices, $\th$ for the $(m+n-1)$th index, and $r$ for the last $(m+n)$th index.  
From the formula in Equation \eqref{eqn:christoffel}, we have the following non-zero entries for the Christoffel symbols 
\begin{align*}
\begin{array}{llll}
    \Gamma^{x_i}_{x_ix_{m-1}} = -\frac{1}{x_{m-1}}, 
    & \Gamma^{x_i}_{x_ir} = \frac{\a_r}{\a},
    & \Gamma^{x_i}_{x_i\theta} = \frac{\a_\th}{\a},
    & \text{for } 1 \leq i \leq m-2. \\[6pt]
    \Gamma^{x_{m-1}}_{x_ix_i} = \frac{1}{x_{m-1}}, 
    & \Gamma^{x_{m-1}}_{x_{m-1}r} = \frac{\a_r}{\a},
    & \Gamma^{x_{m-1}}_{x_{m-1}\theta} = \frac{\a_\th}{\a},
    & \text{for } 1 \leq i \leq m-1. \\[6pt]
    \Gamma^{y_i}_{y_iy_{n-1}} = -\frac{1}{y_{n-1}}, 
    & \Gamma^{y_i}_{y_ir} = \frac{\b_r}{\b},
    & \Gamma^{y_i}_{y_i\theta} = \frac{\b_\th}{\b},
    & \text{for } 1 \leq i \leq n-2. \\[6pt]
    \Gamma^{y_{n-1}}_{y_iy_i} = \frac{1}{y_{n-1}}, 
    & \Gamma^{y_{n-1}}_{y_{n-1}r} = \frac{\b_r}{\b},
    & \Gamma^{y_{n-1}}_{y_{n-1}\theta} = \frac{\b_\th}{\b},
    & \text{for } 1 \leq i \leq n-1. \\[6pt]
    \Gamma^{\theta}_{x_ix_i} = -\frac{\a \a_\th}{x_{m-1}^2f^2}, 
    & \Gamma^{\theta}_{y_iy_i} = -\frac{\b \b_\th}{f^2y_{n-1}^2},
    & \Gamma^{\theta}_{\theta r} = \frac{f'}{f}, 
    & \text{for all } i. \\[6pt]
    \Gamma^{r}_{x_ix_i} = -\frac{\a \a_r}{x_{m-1}^2}, 
    & \Gamma^{r}_{y_iy_i} = -\frac{\b \b_r}{y_{n-1}^2},
    & \Gamma^{r}_{\theta \theta} = -ff',
    & \text{for all } i.
\end{array}
\end{align*}
Recall that the Christoffel symbols are symmetric, that is, $\Gamma^i_{jk}=\Gamma^i_{kj}$. 

By applying the formulas in Equations \eqref{eqn:tensor1} and \eqref{eqn:tensor2}, we obtain the components for the $(4,0)$-curvature tensor. Let us illustrate this process by explicitly working out the computation for $K_\l(X_i, X_j)$ where $1\leq i,j \leq m-1$ and $i\neq j$.  The other tensors can be computed similarly. 
First note that
\begin{align*}
    K_\l(X_i, X_j) &= <R(X_i,X_j)X_i,X_j>= \frac{x^4_{m-1}}{\alpha^4}R_{\scriptstyle{x_ix_jx_ix_j}}=\frac{x^4_{m-1}}{\alpha^4}\sum_l R^l_{\scriptstyle{x_ix_jx_i}}\l_{\scriptstyle{lx_j}}.
\end{align*}

From the matrix in Equation \eqref{eqn:matrix} we have,
\begin{align*}
    \sum_l R^l_{\scriptstyle{x_ix_jx_i}}\l_{\scriptstyle{lx_j}}=R^{x_j}_{\scriptstyle{x_ix_jx_i}}\frac{\alpha^2}{x_{m-1}^2}.
\end{align*}

We apply the formulas in Equations \eqref{eqn:tensor2} and \eqref{eqn:christoffel}, and simplify as the matrix $\Gamma^{x_i}$ has zero entries except $\Gamma^{x_i}_{x_ix_{m-1}}$,$\Gamma^{x_i}_{x_ir}$, and $\Gamma^{x_i}_{x_i\th}$ up to symmetry. 
\begin{align*}
    R^{x_j}_{\scriptstyle{x_ix_jx_i}}&=\sum_l \Gamma^l_{x_ix_i}\Gamma^{x_j}_{x_jl}-\sum_l \Gamma^l_{x_jx_i}\Gamma^{x_j}_{x_il}+\frac{\partial}{\partial x_j}\Gamma^{x_j}_{x_ix_i}-\frac{\partial}{\partial x_i}\Gamma^{x_j}_{x_jx_i}\\
    &=\Gamma^{x_{m-1}}_{x_ix_i}\Gamma^{x_j}_{x_jx_{m-1}}+\Gamma^r_{x_ix_i}\Gamma^{x_j}_{x_jr}+\Gamma^\th_{x_ix_i}\Gamma^{x_j}_{x_j\th}\\
    &=\frac{1}{x_{m-1}}\frac{-1}{x_{m-1}}-\frac{\a \a_r}{x_{m-1}^2} \frac{\a_r}{\a}-\frac{\a\a_\th}{f^2 x_{m-1}^2}\frac{\a_\th}{\a}\\
    &=-\frac{1}{x^2_{m-1}}\left( 1+\a_r^2+\frac{\a_\th^2}{f^2}\right).
\end{align*}

Therefore we have 
\begin{align*}
    K_\l(X_i, X_j) &= \frac{x^4_{m-1}}{\alpha^4}R^{x_j}_{\scriptstyle{x_ix_jx_i}}\frac{\alpha^2}{x^2_{m-1}}\\
    &=-\frac{1}{\a^2}\left(1+\a_r^2+\frac{\a_\th^2}{f^2}\right).
\end{align*}

The second (equivalent) method to calculate the formulas in Theorem \ref{thm: curvature formula polar} is to calculate the Lie brackets for the orthonormal basis in Equations \eqref{eqn:ON basis 1}--\eqref{eqn:ON basis 4}, use these to calculate the Levi-Civita connection, and then use the connection to calculate the components of the sectional curvature tensor. 

The relevant Lie brackets needed for Equation \eqref{eqn:mixed term} are
    \begin{equation*}
    [X_i, Z_\th] = \frac{\a_\th}{\a f} X_i, \quad [X_i, Z_r] = \frac{\a_r}{\a} X_i, \quad [Z_\th, Z_r] = - \frac{f'}{f} Z_\th.
    \end{equation*}
Using the Koszul formula for an orthonormal basis
    \begin{equation*}
    \la \nabla_A B, C \ra = - \frac{1}{2} \left( \la [B,C], A \ra + \la [A,C], B \ra + \la [B,A],C \ra \right),
    \end{equation*}
one calculates
    \begin{equation*}
    \nabla_{X_i} X_i = - \frac{\a_r}{\a} Z_r - \frac{\a_\th}{\a f} Z_\th, \quad \nabla_{Z_r}{X_i} = 0, \quad \nabla_{Z_r} Z_r = 0, \quad \nabla_{Z_r} Z_\th = 0.
    \end{equation*}
Let us quickly note that no Lie brackets other than what are listed above contribute to the above formulas for the connection.  

Finally, we compute
    \begin{align*}
     \la R_\l (X_i, Z_r) X_i, Z_\th \ra_\l &= \la \nabla_{Z_r} \nabla_{X_i} X_i - \nabla_{X_i} \nabla_{Z_r} X_i + \nabla_{[X_i, Z_r]} X_i, Z_{\th} \ra_\l \\
     &= \left\langle \nabla_{\frac{\partial}{\partial r}} \left( - \frac{\a_r}{\a} Z_r - \frac{\a_\th}{\a f} Z_\th \right) - \nabla_{X_i} (0) + \frac{\a_r}{\a} \nabla_{X_i} X_i, Z_\th \right\rangle_{\l} 
     \end{align*}
     \begin{align*}
     &= \left\langle \frac{-\a_{rr} \a + (\a_r)^2}{\a^2} Z_r - \frac{\a_{r \th} f \a - \a_\th (f' \a + f \a_r)}{\a^2 f^2} Z_\th + \frac{\a_r}{\a} \left( - \frac{\a_r}{\a} Z_r - \frac{\a_\th}{\a f} Z_\th \right), Z_\th \right\rangle_{\l} \\
     &= \frac{f'}{\a f^2} \a_\th - \frac{1}{\a f} \a_{r \th} = \frac{1}{\a f} \left( \frac{\a_\th f'}{f} - \a_{r \th} \right).
    \end{align*}

Lastly, note that the above calculation also gives the formula for $K_\l (X_i, Z_r)$.  

\end{proof} 

\begin{cor}\label{cor: only nonzero term}
    Using the notation of Theorem \ref{thm: curvature formula polar}, when $\a = \b = 1$ the only nonzero components of the curvature tensor for $\tilde{\l}$ are
    \begin{equation*}
        K_\l (X_i, X_j) = -1,\quad K_\l (Y_i, Y_j) = -1,\quad K_\l (Z_r, Z_\th) = - \frac{f''}{f}.
    \end{equation*}
\end{cor}

\section{Proof of Theorem \ref{thm:main theorem 2}}\label{sect:proof of main theorem 2}

\subsection{Turning $S$ into a cusp of $N$}
To prove Theorem \ref{thm:main theorem 2} $(1)$, we need to warp the product hyperbolic metric near $S$ in order to turn $S$ into a cusp of $N$.
We do this by defining functions $\a$ and $\b$ on $(- \infty, \eta) \times \S^1$ and $f$ on $(- \infty, \eta)$ in Equation \eqref{eqn:polar coordinate warped product metric} for some sufficiently small $\eta > 0$.  
As the radial parameter $r$ approaches $- \infty$, this has the geometric effect of turning $S$ into a cusp of $N$. 

\begin{remark}\label{rmk:finite-volume}
    By Fubini's Theorem, and since $\S^1$ is compact, a sufficient condition for the metric $\tilde{\l}$ in Equation \eqref{eqn:polar coordinate warped product metric} to descend to a metric on $N = M \setminus S$ with finite-volume is that all three functions $\a, \b,$ and $f$ are bounded for $r \in (- \infty, \e)$, and at least one of the functions decays to $0$ exponentially fast as $r \to - \infty$. 
\end{remark}

Let $E = \H^{m-1} \times \H^{n-1} \times \S^1$, and consider the metric $\tilde{\l}$ from Equation \eqref{eqn:polar coordinate warped product metric} defined on $E \times (0, \infty)$.
Let $q \in E \times (0, \infty)$, let $\sigma \subset T_q (E \times (0, \infty))$ be a $2$-plane, and consider the orthonormal frame about $q$ given by Equations \eqref{eqn:ON basis 1}--\eqref{eqn:ON basis 4}.
Choose an orthonormal basis $(A, B)$ for $\s$ which satisfies that $A$ is orthogonal to $Z_r$, which is always possible since the orthogonal complement to $Z_r$ has codimension one.

Choose a unit vector $X_1$ in the $\H^{m-1}$ component of $T_q (E \times (0, \infty))$ and a unit vector $Y_1$ in the $\H^{n-1}$ component of $T_q (E \times (0, \infty))$ which satisfy that $X_1$ is parallel to $\text{proj}_{\H^{m-1}}(A)$ and $Y_1$ is parallel to $\text{proj}_{\H^{n-1}}(A)$.
Finally, choose corresponding vectors $X_2$ and $Y_2$ such that $\text{proj}_{\H^{m-1}}(\s) \subseteq \text{span}(X_1, X_2)$ and $\text{proj}_{\H^{n-1}}(\s) \subseteq \text{span}(Y_1, Y_2)$.
Then there exist constants $a_1, a_3, a_5, b_1, b_2, \hdots , b_6$ such that
\begin{equation*}
A = a_1 X_1 + a_3 Y_1 + a_5 Z_\th \quad \text{and} \quad B = b_1 X_1 + b_2 X_2 + b_3 Y_1 + b_4 Y_2 + b_5 Z_\th + b_6 Z_r.
\end{equation*}
Note that, since $A$ and $B$ are orthonormal, we have
\begin{equation*}
a_1^2 + a_3^2 + a_5^2 = 1, \qquad b_1^2 + \hdots + b_6^2 = 1, \qquad a_1 b_1 + a_3 b_3 + a_5 b_5 = 0.
\end{equation*}

In the sequel we use the following notation.  Let
\begin{equation*}
\G_1 = X_1, \quad \G_2 = X_2, \quad \G_3 = Y_1, \quad \G_4 = Y_2, \quad \G_5 = Z_\th, \quad \G_6 = Z_r,
\end{equation*}
and
\begin{equation*}
R_{ijk \ell} = \la R(\G_i, \G_j)\G_k, \G_\ell \ra_{\l}.
\end{equation*}
We then compute
\begin{align}
K_\l(\s) &= \la R(A, B)A, B \ra_{\l} = a_1^2 b_2^2 R_{1212} + (a_1 b_3 - a_3 b_1)^2 R_{1313} + a_1^2 b_4^2 R_{1414} \notag \\
&+ a_3^2 b_2^2 R_{2323} + a_5^2 b_2^2 R_{2525} + a_3^2 b_4^2 R_{3434} + a_5^2 b_4^2 R_{4545} + a_5^2 b_6^2 R_{5656} \notag \\
&+ (a_1 b_5 - a_5 b_1)^2 R_{1515} + 2a_1 b_6 (a_1 b_5 - a_5 b_1) R_{1516} + a_1^2 b_6^2 R_{1616} \label{eqn:quad form 1} \\
&+ (a_3 b_5 - a_5 b_3)^2 R_{3535} + 2a_3 b_6 (a_3 b_5 - a_5 b_3) R_{3536} + a_3^2 b_6^2 R_{3636} \label{eqn:quad form 2}.
\end{align}
Lines \eqref{eqn:quad form 1} and \eqref{eqn:quad form 2} can be rewritten as matrix products
\begin{equation}\label{eqn:matrix equation 1}
\begin{bmatrix}
a_1 b_5 - a_5 b_1 & a_1 b_6
\end{bmatrix}
\begin{bmatrix}
R_{1515} & R_{1516} \\
R_{1516} & R_{1616}
\end{bmatrix}
\begin{bmatrix}
a_1 b_5 - a_5 b_1 \\
a_1 b_6
\end{bmatrix},
\end{equation}
and
\begin{equation}\label{eqn:matrix equation 2}
\begin{bmatrix}
a_3 b_5 - a_5 b_3 & a_3 b_6
\end{bmatrix}
\begin{bmatrix}
R_{3535} & R_{3536} \\
R_{3536} & R_{3636}
\end{bmatrix}
\begin{bmatrix}
a_3 b_5 - a_5 b_3 \\
a_3 b_6
\end{bmatrix}.
\end{equation}

With these calculations we are ready to prove the following proposition, which gives a relatively simple criterion to determine exactly when $\tilde{\l}$ is non-positively curved.

\begin{proposition}\label{prop:curvature conditions}
The metric $\tilde{\l}$ in Equation \eqref{eqn:polar coordinate warped product metric} is non-positively curved if and only if both of the following two conditions hold:
	\begin{enumerate}
	\item All of the curvatures in Equations \eqref{eqn:curvature component}--\eqref{eqn:curvature component 2} are non-positive.

	\item All matrices of the form
	\begin{equation}\label{eqn:matrices 1}
	M_1:=\begin{bmatrix}
	K_\l (X_i, Z_\th) & <R_\l(X_i,Z_r)X_i,Z_\th>_\l \\
	<R_\l(X_i,Z_r)X_i,Z_\th>_\l & K_\l(X_i, Z_r)
	\end{bmatrix}
	\end{equation}
	and
	\begin{equation}\label{eqn:matrices 2}
	M_2:=\begin{bmatrix}
	K_\l (Y_i, Z_\th) & <R_\l(Y_i,Z_r)Y_i,Z_\th>_\l \\
	<R_\l(Y_i,Z_r)Y_i,Z_\th>_\l & K_\l(Y_i, Z_r)
	\end{bmatrix}
	\end{equation}
are negative semi-definite.
	\end{enumerate}
\end{proposition}

\begin{proof}  
First, notice that the matrices in Equations \eqref{eqn:matrices 1} and \eqref{eqn:matrices 2} are of the exact same form as the matrices in Equations \eqref{eqn:matrix equation 1} and \eqref{eqn:matrix equation 2}.  

To prove the proposition, we first assume that conditions $(1)$ and $(2)$ hold.  
Then, by $(2)$, lines \eqref{eqn:quad form 1} and \eqref{eqn:quad form 2} in the calculation for $\la R(A,B)A, B \ra_\l$ are non-positive.  
By $(1)$, the first two lines are also non-positive.  
Therefore, the metric $\tilde{\l}$ must have non-positive sectional curvature.

So let us now assume that all sectional curvatures of $\tilde{\l}$ are non-positive.  
This clearly implies statement (1).   
To verify statement $(2)$, we consider specific $2$-planes $\s$ which correspond to choosing specific coefficients $a_1, a_3, a_5, b_1, b_2, \hdots b_6$.  
In particular, we begin by specifying that
\begin{equation*}
a_3 = a_5 = b_1 = b_2 = b_3 = b_4 = 0.
\end{equation*}
The only nonzero line in the calculation of $\la R(A,B)A, B \ra_\l$ with this choice of coefficients is \eqref{eqn:quad form 1}.  
Note that $a_3 = a_5 = 0$ forces $a_1 = \pm 1$.
The vector
\begin{equation*}
\begin{bmatrix}
a_1 b_5 \\
a_1 b_6
\end{bmatrix}
= 
\pm \begin{bmatrix}
b_5 \\ b_6
\end{bmatrix}
\end{equation*}
can realize any direction in $\R \mathbb{P}^1$ for different choices of $b_5$ and $b_6$.  
So the only way to ensure that \eqref{eqn:quad form 1} is non-positive for all possible selections of these variables is if the matrix
\begin{equation*}
M_1=\begin{bmatrix}
R_{1515} & R_{1516} \\
R_{1516} & R_{1616}
\end{bmatrix}
\end{equation*}
is negative semi-definite.  
An identical analysis shows that we also need the matrix
\begin{equation*}
M_2=\begin{bmatrix}
R_{3535} & R_{3536} \\
R_{3536} & R_{3636}
\end{bmatrix}
\end{equation*}
to be negative semi-definite, verifying the necessity of condition $(2)$.  
\end{proof}

A similar argument can be used to provide a criterion guaranteeing a small positive upper bound on the sectional curvatures.

\begin{proposition}\label{cor: almost nonpositive curvature}
    Assume that there are universal constants $\kappa_{max}\geq 0$ and $ \mu_{max}\geq 0$, such that the metric $\tilde{\lambda}$ in Equation \eqref{eqn:polar coordinate warped product metric} has the following two properties:
    \begin{enumerate}
        \item all of the curvatures in Equations \eqref{eqn:curvature component}--\eqref{eqn:curvature component 2} are $\leq \kappa_{max}$, and
        \item the largest eigenvalues of the matrices $M_1, M_2$ are $\leq \mu_{max}$.
    \end{enumerate}
    Then the metric $\tilde{\l}$ has all sectional curvatures $\leq 11\kappa_{max} + 10\mu_{max}$.
\end{proposition}

\begin{proof}
    Looking at the expression for the sectional curvatures, and remembering that all the coefficients $a_i, b_i$ have norm $\leq 1$, we obtain an upper bound of the form:
    \begin{align}
        K_\lambda(\sigma) &\leq 11\kappa_{max} +v_1M_1v_1^T + v_2M_2v_2^T \notag, \\
        &\leq 11\kappa_{max} + ||v_1||^2\mu(M_1)+ ||v_2||^2\mu(M_2) \label{eqn:eigval-upper-bd}.
    \end{align}
    Here $v_1=\begin{bmatrix} a_1 b_5 - a_5 b_1 & a_1 b_6 \end{bmatrix}$, $v_2=\begin{bmatrix} a_3 b_5 - a_5 b_3 & a_3 b_6\end{bmatrix}$, and $\mu(M_i)$ denotes the largest eigenvalue of the matrices $M_i$. By hypothesis, we have $\mu(M_i)\leq \mu_{max}$. Since the $a_i, b_i$ lie in the interval $[-1,1]$, the vectors $v_1, v_2$ have norm squared bounded by $||v_i||^2\leq 5$, establishing the Proposition.
\end{proof}

In other words, Proposition \ref{cor: almost nonpositive curvature} shows that to bound the sectional curvatures above by a small positive constant, it suffices to (uniformly) bound both (1) the curvatures on the coordinate hyperplanes and (2) the largest eigenvalues of the matrices $M_1, M_2$, by a small positive constant. We are now ready to prove Theorem \ref{thm:main theorem 2} (1).

\begin{proof}[Proof of Theorem \ref{thm:main theorem 2} (1)]

We again let $E = \H^{m-1} \times \H^{n-1} \times \S^1$, and let $\e > 0$ be given.
Choose $\eta \in (0,1)$ to be less than the normal injectivity radius of $S$ in $M$.
We first define a piecewise-smooth metric $\tilde{g}'$ by prescribing specific values for the functions $\a$, $\b$, and $f$ in Equation \eqref{eqn:polar coordinate warped product metric} on different portions of $E \times \R$.  
We will show that, for a sufficiently small positive parameter $\d < \eta$, the sectional curvatures of $\tilde{g}'$ at all points where it is smooth are bounded above by $\e/2$.  
We will then argue that we can $C^2$-approximate $\tilde{g}'$ by a smooth metric $\tilde{g}$ whose sectional curvatures increase by at most $\e/2$ (but whose lower curvature bound decreases without bound as $\e \to 0$, forcing $\d \to 0$).  
The metric $\tilde{g}$ will agree with $\tilde{\h}_m \oplus \tilde{\h}_n$ outside of a sufficiently small tubular neighborhood of $\H^{m-1} \times \H^{n-1}$, and is invariant under isometries of $\tilde{\h}_m \oplus \tilde{\h}_n$ that stabilize $\H^{m-1} \times \H^{n-1}$. Thus $\tilde g$, restricted to the small tubular neighborhood of $\H^{m-1} \times \H^{n-1}$, will descend 
to a well-defined metric $g$ on a (punctured) small tubular neighborhood of $S$ in $M$. Since the metric $g$ coincides with the locally symmetric metric near the boundary of the tubular neighborhood, we can replace the locally symmetric metric on the tubular neighborhood by the metric $g$. This has the effect of puncturing the manifold by removing $S$, and inserting a cusp supported in the tubular neighborhood of $S$.

\medspace

As mentioned above, the metric $\tilde{g}'$ will depend on a small parameter $\d > 0$. 
Define functions $f , \v \colon \R \to \R$ by
\begin{equation*}
f(r)=\begin{cases}
r &r \ge \eta\\
r &\d<r<\eta\\
\d &0<r\le \d\\
\d e^{\d r} &r \le 0
\end{cases},
\qquad\qquad
\v(r)=\begin{cases}
1 & r \ge \eta\\
h(r)&\d <r < \eta\\
0 & 0<r \le \d\\
0 &r \le 0
\end{cases},
\end{equation*}
where 
\begin{equation*}
h(r) = \frac{1}{\ln \left( \frac{\eta}{\d} \right) } \ln \left( \frac{r}{\d} \right).
\end{equation*}
With the above definition of $h(r)$, one immediately sees that $\v$ is continuous at both $r = \d$ and $r = \eta$. 
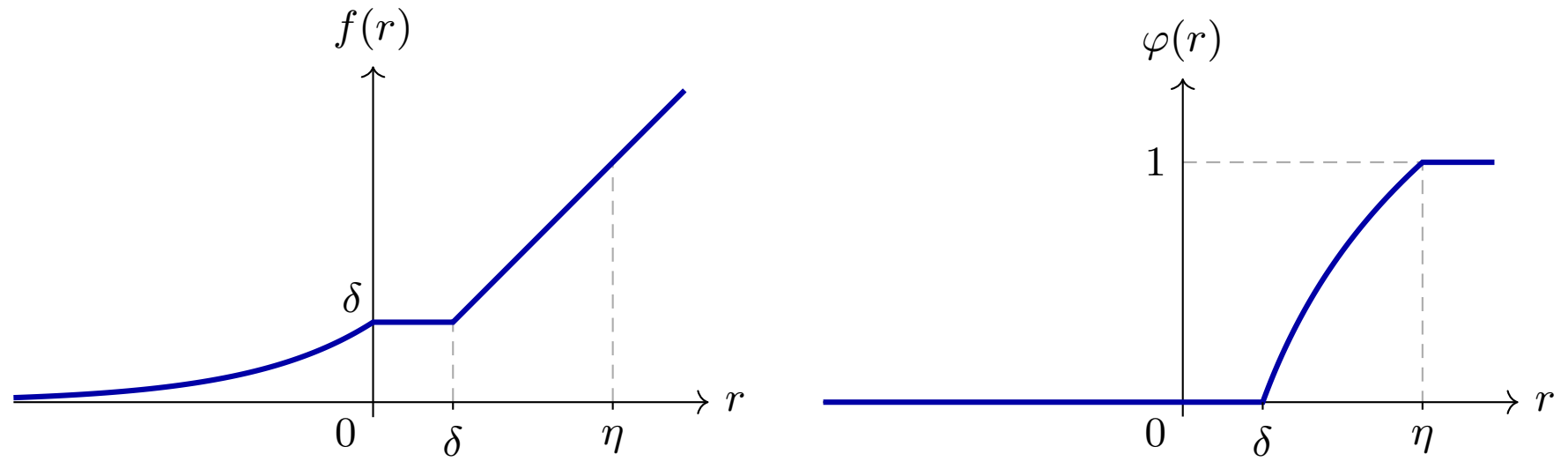
\begin{figure}
    \centering

\begin{tikzpicture}[
scale=0.8,
  graph/.style={draw=blue!65!black, line width=1.15pt},
  axis/.style={->, thin},
  guide/.style={densely dashed, gray!65, thin},
  every node/.style={font=\small}
]
 
  \def\Dval{0.8}
  \def\Eval{2.4}
  \pgfmathsetmacro{\rMin}{-1.5*\Eval}
  \pgfmathsetmacro{\rMax}{1.3*\Eval}
  \pgfmathsetmacro{\rEnd}{1.4*\Eval}
  \pgfmathsetmacro{\unitlen}{2.25/\Eval}

  \begin{scope}[x=\unitlen cm, y=\unitlen cm]
    \draw[guide] (\Dval,0) -- (\Dval,\Dval);
    \draw[guide] (\Eval,0) -- (\Eval,\Eval);
    \draw[axis] (\rMin,0) -- (\rEnd,0) node[right] {$r$};
    \draw[axis] (0,-0.15) -- (0,\rEnd) node[above] {$f(r)$};

    \draw[graph] plot[domain=\rMin:0,samples=100,variable=\t]
      (\t,{\Dval*exp(\Dval*\t)})
      -- (\Dval,\Dval) -- (\rMax,\rMax);

    \node[below left] at (0,0) {$0$};
    \node[above left, inner sep=2pt] at (0,\Dval) {$\delta$};
    \foreach \r/\lab in {\Dval/\delta,\Eval/\eta}
      \draw (\r,0) -- ++(0,-2pt) node[below] {$\lab$};
  \end{scope}

  \begin{scope}[xshift=7.6cm, x=\unitlen cm, y=2.25cm]
    \draw[guide] (0,1) -- (\Eval,1) -- (\Eval,0);
    \draw[axis] (\rMin,0) -- (\rEnd,0) node[right] {$r$};
    \draw[axis] (0,-0.06) -- (0,1.35) node[above] {$\varphi(r)$};

    \draw[graph] (\rMin,0) -- (\Dval,0)
      -- plot[domain=\Dval:\Eval,samples=100,variable=\t]
         (\t,{ln(\t/\Dval)/ln(\Eval/\Dval)})
      -- (\rMax,1);

    \node[below left] at (0,0) {$0$};
    \node[left] at (0,1) {$1$};
    \foreach \r/\lab in {\Dval/\delta,\Eval/\eta}
      \draw (\r,0) -- ++(0,-2pt) node[below] {$\lab$};
  \end{scope}
\end{tikzpicture}
    \caption{Graphs of the functions $f$ and $\phi$}
    \label{fig:placeholder}
\end{figure}

For future calculations, quickly note that
\begin{equation*}
f'(r)=\begin{cases}
1 &r > \eta\\
1 &\d<r<\eta\\
0 &0<r < \d\\
\d^2 e^{\d r} &r < 0
\end{cases},
\qquad\qquad
\v'(r)=\begin{cases}
0 & r > \eta\\
h'(r)&\d <r < \eta\\
0 & 0<r < \d\\
0 &r < 0
\end{cases},
\end{equation*}
\begin{equation*}
h'(r) = \frac{1}{r \ln \left( \frac{\eta}{\d} \right) }, \qquad \qquad h''(r) =  \frac{-1}{r^2 \ln \left( \frac{\eta}{\delta} \right)} .
\end{equation*}
Also notice that, as $\d \to 0$, $f$ is approximately $C^2$ at $r = 0$ and $\v$ is approximately $C^2$ at $r = \eta$.
However, $f$ and $\v$ are both far from being $C^1$ at $r = \d$.  

We define the metric $\tilde{g}'$ via Equation \eqref{eqn:polar coordinate warped product metric}.
We of course use the function $f(r)$ above for the coefficient of $d \th^2$.  
For $\a$ and $\b$ we define
\begin{equation*}
\a(r, \th) = \cosh(\v(r) f(r) \cos(\th)), \qquad \qquad \b(r, \th) = \cosh(\v(r) f(r) \sin(\th)).
\end{equation*}
Note that $\tilde{g}' = \h_m \oplus \h_n$ for $r > \eta$. We now verify that the singular metric $\tilde{g}'$, with parameter $\d$ chosen small enough, yields a metric with all sectional curvatures $<\e/2$.

\medskip

\noindent\underline{\bf Claim 1}: Consider the singular metric $\tilde{g}'$ away from the singular values $r = 0 , \d ,$ and $\eta$. Then when $\delta$ is small enough, the singular metric has all coordinate plane curvatures described in Equations \eqref{eqn:curvature component}--\eqref{eqn:curvature component 2} that are $\leq \e/44$. 

\medskip

Of the sectional curvatures in Equations \eqref{eqn:curvature component}--\eqref{eqn:curvature component 2}, it is clear that \eqref{eqn:curvature component}, \eqref{eqn:curvature component 3}, and \eqref{eqn:curvature component 2} are all non-positive. Indeed, the expressions in \eqref{eqn:curvature component} and \eqref{eqn:curvature component 3} are
$\leq 0$ regardless of the choice of the functions $\v$ and $f$, while the sign of the expression in \eqref{eqn:curvature component 2} depends only on the concavity of $f$ (and our function satisfies $f''\geq 0$ away from the singular sets). 

For the sectional curvature in \eqref{eqn:curvature component 8}, an easy calculation gives us
$$K(X_i,Y_j)=\Big(\cos(\theta)\tanh(\v f \cos(\theta))\Big)\Big(\sin(\theta)\tanh(\v f \sin(\theta))\Big)\Big(\v^2 -(\v 'f + \v f')^2\Big).$$
The first two terms in this product are always $\geq 0$, so the sign of this expression only depends on the third term. But this third term evaluates to zero, except when $r\in (\d, \eta)$ where it equals
$$\v^2 -(\v 'f + \v f')^2=
-2rhh'-r^2(h')^2.$$
Since $h, h'\geq 0$ on the interval $ (\d, \eta)$, we conclude that the curvature in \eqref{eqn:curvature component 8} is indeed $\leq 0$.

\medskip

More work is required to estimate the curvature terms in Equations \eqref{eqn:curvature component 4} -- \eqref{eqn:curvature component 7}. These are
the four terms
$K_{g'}(X_i, Z_\th)$, $K_{g'}(X_i, Z_r)$, $K_{g'}(Y_i, Z_{\th})$, $K_{g'}(Y_i, Z_r)$. For the remainder of the proof we only consider calculations with the $X_i$, as the $Y_i$ calculations are identical.

A direct calculation using the formulas in Theorem \ref{thm: curvature formula polar} shows that the two terms in Equations \eqref{eqn:curvature component 4} and \eqref{eqn:curvature component 6} reduce to:
\begin{align*}
& K_{g'}(X_i,Z_\theta)
  = -\v^2\sin^2(\theta)
     -\boxed{\left(\v'f'+\v\frac{(f')^2-1}{f}\right)\cos(\theta)\tanh(\v f\cos(\theta))}, \\
& K_{g'}(X_i,Z_r)
  = -(\v'f+\v f')^2\cos^2(\theta)
     -\boxed{{(\v''f+2\v'f'+\v f'')\cos(\theta)\tanh(\v f\cos(\theta))}}.\\
\end{align*}
Note that the unboxed portions of the two equations above are clearly $\leq 0$. Thus the two boxed terms are the {\it error terms}, the only terms that could potentially lead to one of these curvature expressions being positive. Moreover, from the choice of our functions $f$, $\v$, the boxed terms are non-zero only when $r\in (\d, \eta)$. This already tells us that the sectional curvatures in Equations \eqref{eqn:curvature component 4} and \eqref{eqn:curvature component 6} are non-positive, except possibly in the region $ (\d, \eta)$. 

To control the curvature on this interval, we calculate an upper bound on the error terms as follows. First note that from our choice of functions we have
$$\v((f')^2-1)\equiv0,$$ simplifying the first error term. Next recalling that $| \tanh(x)| \leq |x|$ for all $x$, gives us the following inequalities for the two error terms
\begin{align*}
\left|\left(\v'f'\right)\cos(\theta)\tanh(\v f\cos(\theta))\right| &\le |\v'||f'||\v||f|,\\
&=|rh(r)h'(r)|,\\
\left|(\v''f+2\v'f'+\v f'')\cos(\theta)\tanh(\v f\cos(\theta))\right|&\le |\v''f+2\v'f'+\v f''||\v||f| ,\\
&=|rh(r)||rh''(r)+2h'(r)|,
\end{align*}
where for the equalities, we have used that $f(r)=r$, $\v(r)=h(r)$ on the interval $(\d, \eta)$. We can now use the specific form of our function $h(r)$ to give upper bounds:

\begin{equation}\label{eqn:bounding inequality 1}
|rh(r)h'(r)|\le \left|\ \frac{1}{\ln(\eta/\d)}\cdot \ln(r/\d)\right |\left|\frac{1}{\ln(\eta/\d)}\right|\le  \left|\frac{1}{\ln(\eta/\d)}\right| ,
\end{equation}
\begin{equation}\label{eqn:bounding inequality 2}
|rh(r)||rh''(r)+2h'(r)|\le \left|\frac{r}{\ln(\eta/\d)}\cdot\ln(r/\d)\right |\left|\frac{1}{r\cdot\ln(\eta/\d)}\right| \le\left|\frac{1}{\ln(\eta/\d)}\right|.
\end{equation}
But as $\d\to0^+$, the term $\frac{1}{\ln(\eta/\d)}$ tends to zero. It follows that by choosing $\d > 0$ sufficiently small, the curvatures in Equations \eqref{eqn:curvature component 4} and \eqref{eqn:curvature component 6} can be bounded above by $\e/44$, as desired. As mentioned earlier, the estimates for the sectional curvatures in \eqref{eqn:curvature component 5} and \eqref{eqn:curvature component 7} are identical. This completes the proof of {\bf Claim 1}.

\medskip

\noindent \underline{\bf Claim 2:} For $\delta$ sufficiently small, the singular metric $\tilde{g}'$, away from the singular sets $r=0,\delta, \eta$, has the property that the eigenvalues of the corresponding matrices $M_1, M_2$ are all $\leq \e /40$.

\medskip

We will verify this for the matrix $M:=M_1$, the computations for $M_2$ being completely analogous. Recall that $M$ is described in Equation \eqref{eqn:matrices 1}. The diagonal entries are given by the curvature terms in \eqref{eqn:curvature component 4} and \eqref{eqn:curvature component 6}, which were analyzed earlier. The off-diagonal entries are given by the mixed term described in \eqref{eqn:mixed term}. For our metric $\tilde{g}'$, the mixed term takes the form
$$\langle R_{g'}(X_i,Z_r)X_i,Z_\theta\rangle = \v(\v'f+\v f')\sin(\theta)\cos(\theta)+\boxed{{\v'\sin(\theta)\tanh(\v f\cos(\theta))}},$$
where as before we have a boxed {\it error term}. 
As in the earlier cases, the boxed error term vanishes for the locally symmetric metric.
Note that $\v'$ vanishes except when $r\in (\d, \eta)$. Within this range, the boxed error term can be estimated above by
$$|\v'\sin(\theta)\tanh(\v f\cos(\theta))| \le |\v '||\v||f|=|rh(r)h'(r)|.$$

It follows from the estimate in Equation \eqref{eqn:bounding inequality 1} that this upper bound can be made arbitrarily small by choosing $\delta$ sufficiently close to zero. 

We now consider the matrices $M$ over each of the four regions of interest.  Substituting in our specific functions $f, \v$, the matrices are given by:
\begin{align*}
    & M = {\begin{bmatrix}
-\sin^2(\theta) & \sin(\theta)\cos(\theta)\\
\sin(\theta)\cos(\theta) & -\cos^2(\theta)
\end{bmatrix}} & \eta < r \\
    & M \approx M'= {\begin{bmatrix}
-h^2\sin^2(\theta) & h(rh'+h)\sin(\theta)\cos(\theta)\\
h(rh'+h)\sin(\theta)\cos(\theta) & -(rh'+h)^2\cos^2(\theta)
\end{bmatrix}} & \d<r<\eta \\
    & M = {\begin{bmatrix}
0 & 0\\
0 & 0
\end{bmatrix}} & r<\d, r\neq 0
\end{align*}

\noindent where for the interval $(\d, \eta)$ we have only given an {\it approximation} $M'$ of $M$ where we have dropped out the error term. When $r\in (0, \delta)$ and $r\in (-\infty, 0)$, $M$ is the zero matrix, so has largest eigenvalue $=0$. We now observe that on the two intervals $ (\delta, \eta)$ and $ (\eta, \infty)$, these matrices have non-positive trace and zero determinant. This immediately implies that these matrices also have largest eigenvalue $=0$. 

\medskip

This gives the desired eigenvalue bound, except potentially for the region when $r\in (\d, \eta)$, where the matrices $M'$ above are only approximations for $M$. The actual entries in the matrix $M$ should include the boxed terms, which we have seen perturbs all the entries of $M'$ by some amount that is $O\left(\frac{-1}{\ln(\d)}\right)$. The largest eigenvalue of a $2\times 2$ matrix varies continuously as a function of the matrix entries. So for any single matrix, we can choose $\delta$ small enough to keep the largest eigenvalue $\leq \e /40$. But $M'$ is a continuous family of matrices, parametrized by $r\in (\d, \eta)$, $\theta\in [0,2\pi]$, and we need to ensure that we can pick a uniform $\delta$ for the entire family\footnote{The functions $f$ and $\v$ are continuous over the compact interval $[\d, \eta]$, but $\v' \to \infty$ as $\d \to 0$. So it is a priori possible that the $L^1$-estimate gets worse as $\d \to 0$. This is why further analysis is required here. }  .

To this end, recall that the characteristic polynomial of a $2 \times 2$ matrix $M$ is 
$$\chi(x)=x^2 - \text{tr}(M)x + \det(M),$$
and hence the largest root will be bounded by
\begin{align*}
    \mu(M)&=\frac{1}{2} \left( \text{tr}(M) +\sqrt{\text{tr}(M)^2-4\det(M)} \right),\\
    &\leq \text{tr}(M)+|\text{tr}(M)|+2\sqrt{|\det(M)|} .
\end{align*}
If we perturb $M'$ to a new matrix $M$ by changing the entries by some small amount $O(\epsilon)$, the two trace terms change by $O(\epsilon)$ (with uniform constant). From the formula for the determinant, we see that $\det(M)=\det(M')+O_{M'}(\epsilon)$, where the constant in the error term depends on the entry-wise $L^1$-norm of the matrix $M'$, i.e. the sum of the absolute value of the matrix entries. When $r\in (\d, \eta)$, we have $h(r)\in (0,1)$ and 
$$rh'(r)+h(r) = \frac{1+\ln(r/\d)}{\ln(\eta/\d)} \leq 2,$$
so all the matrix entries in our approximation $M'$ have norm $\leq 2$. Since the approximate matrices have $\det(M')=0$, we conclude that the actual matrices $M$ have $\det(M) = O(\epsilon)$, with uniform constant, on the domain $ (\d, \eta)$. Finally, we can substitute in these bounds, recalling that the approximate matrices satisfy $\text{tr}(M')\leq 0$ and $\det(M')=0$, to get that the actual matrices satisfy the estimate:
$$\mu_{max}(M) = O(\epsilon) + O(\sqrt{\epsilon}) = O\left(\sqrt{\frac{-1}{\ln{\d}}}\right),$$
where the second equality follows from the fact that $\epsilon=O\left(\frac{-1}{\ln{\d}}\right)$, and since $0<\epsilon<1$ we also have $\epsilon<\sqrt{\epsilon}$. As this estimate holds uniformly on the region where $r\in (\d, \eta)$, taking $\d$ small enough  guarantees the largest eigenvalue of the matrices $M$ are $\leq \e/40$, completing the proof of {\bf Claim 2}.

\medskip 

We can now set $\kappa_{max}:= \e/44$ and $\mu_{max}:= \e /40$.
Then the estimates obtained in {\bf Claim 1} and {\bf Claim 2} show that the hypotheses of Proposition \ref{cor: almost nonpositive curvature} are satisfied. Hence for $\delta$ small enough, {\bf all} sectional curvatures for the singular metric $\tilde{g}'$ are $\leq \epsilon/2$. The last step is then to smooth the singularities while not introducing too much additional positive curvature.

\begin{remark} \label{rmk:shrinking delta}
    Observe that, as $\delta\rightarrow 0^+$, all the upper bounds in {\bf Claim 1} and {\bf Claim 2} get closer to zero. This means that in subsequent arguments, we still have the freedom of further shrinking $\delta$ if needed.
\end{remark}

\medskip

\noindent \underline{\bf Claim 3:} The singular metric $\tilde{g}'$ can be smoothed near the three singular sets $r=0, \delta, \eta$, to a smooth metric $\tilde g$. The smoothing can be chosen to still satisfy the hypotheses of Proposition \ref{cor: almost nonpositive curvature}, with $\kappa_{max}:= \e/22$ and $\mu_{max}:= \e /20$. Thus the smoothed metric $\tilde{g}$ has all curvatures $\leq \e$.

\medskip

In order to smooth the singular metric $\tilde{g}'$, we just need to smooth the function $f$ near $r = 0$ and $\d$, and smooth $\v$ near $r=\d$ and $\eta$. In order to control $\kappa_{max}$ and $\mu_{max}$, we will need to analyze how the perturbations in $f$, $\v$ affect the expressions for (1) the sectional curvature along the coordinate hyperplanes, and (2) the largest eigenvalues of the matrices $M_i$. We analyze each singular value in turn.

\medskip

\noindent {\underline{Singularity at $r=0$}:} In a neighborhood of the $r=0$ singularity, note that $\v \equiv 0$ (so is smooth). Thus we only need to smooth the function $f$. Near the singularity at $r=0$, the function $f$ is given by:
$$f(r)= \begin{cases}
    \delta & r\geq 0\\
    \delta e^{\delta r} & r\leq 0
\end{cases}
$$
and hence is $C^0$, but not $C^1$ at the singularity. Any smoothing $f_\epsilon$ will $C^0$-approximate $f$, so in particular $f_\epsilon$ will be bounded close to $\delta$. In terms of $f_\epsilon '$, this is expected to transition from $\approx\delta^2$ to zero. 

Considering the expressions for the coordinate plane curvatures, we note that Equations (\ref{eqn:curvature component}) and (\ref{eqn:curvature component 3}) are $\leq 0$ regardless of $f$, $\v$. The vanishing of $\v$ also implies that the expressions in Equations (\ref{eqn:curvature component 4})--(\ref{eqn:curvature component 8}) are all $\equiv0$ in a neighborhood of $r=0$, regardless of the choice of smoothing for $f$. This leaves us with Equation (\ref{eqn:curvature component 2}), which we recall states that $K(Z_r, Z_\theta)=-f''/f$. 

Concerning the matrices $M_1$, $M_2$ near the singularity $r=0$, observe that the vanishing $\v \equiv 0$ immediately implies that all the entries vanish (regardless of the smoothing of $f$). Thus the matrices have largest eigenvalue $0$, as desired.

So in order to apply Proposition \ref{cor: almost nonpositive curvature} for the smoothed out metric near $r=0$, it suffices to ensure our local smoothing $f_\epsilon$ near $r=0$ satisfies $-f_\epsilon''/f_\epsilon \leq \e /22$. Such a smoothing is provided by Lemma \ref{lem:smoothing at r=0} in the Appendix.

\medskip 

\medskip

\noindent{\underline{Singularity at $r=\eta$}:} In a neighborhood of the $r=\eta$ singularity, the function $f(r)=r$ is already smooth. So we only need to smooth the function $\v$. Near the singularity at $r=\eta$, the function $\v$ is given by:
$$\v(r)= \begin{cases}
    1 & r\geq \eta\\
    h(r) & r\leq \eta
\end{cases}
$$
and hence is $C^0$, but not $C^1$ at the singularity. A smoothing $\v_\epsilon$ will $C^0$-approximate $\v$, so will be bounded close to $1$ near $\eta$. Moreover, $\v_\epsilon'$ is expected to transition from a value $\approx \frac{1}{\eta\ln{(\eta/\delta)}}$ to zero.

For the coordinate plane curvatures, Equations (\ref{eqn:curvature component}) and (\ref{eqn:curvature component 3}) are $\leq 0$ regardless of $f$, $\v$. The curvature expression in Equation (\ref{eqn:curvature component 2}) is linear in $f''\equiv 0$, hence is also $\equiv 0$. Using the specific form of the function $f$, one can easily check that the remaining coordinate plane curvatures are $\leq \e /22$ provided the smoothing $\v_\epsilon$ 
satisfies the following constraints near $r=\eta$:
\begin{itemize}
    \item $\v_\epsilon' \geq -\e /22$, for control of curvatures in Equations (\ref{eqn:curvature component 4}), (\ref{eqn:curvature component 5}), (\ref{eqn:curvature component 8}),
    \item $r\v_\epsilon ''+ 2\v_\epsilon' \geq -\e /22$, for control of the curvatures in Equations (\ref{eqn:curvature component 6}), (\ref{eqn:curvature component 7}).
\end{itemize}

Observe that if the two expressions above are $\geq 0$, then the sectional curvatures along the coordinate planes would be non-positive. However, an analysis of the function $r^2 \ve'(r)$ shows that there does {\bf not} exist a function that interpolates between $h(r)$ and $1$ and has both these terms $\geq 0$. In particular, any smoothing of $\varphi$ near $r= \eta$ will result in some coordinate hyperplanes having positive curvature near $r= \eta$.

Next we consider the eigenvalues of the matrices $M_i$. As before, we decompose the matrix $M:=M_1$ as a sum $M=M' + \boxed{M}$ by separating the unboxed and boxed portions of each matrix entry:
\begin{align*}
     M' &= {\begin{bmatrix}
-\v_\epsilon^2\sin^2(\theta) & \v_\epsilon(r\v_\epsilon'+\v_\epsilon)\sin(\theta)\cos(\theta)\\
\v_\epsilon(r\v_\epsilon'+\v_\epsilon)\sin(\theta)\cos(\theta) & -(r\v_\epsilon'+\v_\epsilon)^2\cos^2(\theta)
\end{bmatrix}}\\
    \boxed{M} &= {\begin{bmatrix}
-(\v_\epsilon')\cos{(\theta)}\tanh{(r\v_\epsilon \cos{\theta})} & {\v_\epsilon'\sin(\theta)\tanh(r\v_\epsilon \cos\theta)}\\
{\v_\epsilon'\sin(\theta)\tanh(r\v_\epsilon \cos\theta)} & {-(r\v_\epsilon''+2\v_\epsilon')\cos(\theta)\tanh(r\v_\epsilon \cos\theta)}
\end{bmatrix}}\\
\end{align*}
The approximation $M'$ has determinant zero, and non-positive trace, so has largest eigenvalue $=0$. As discussed in {\bf Claim 2}, the largest eigenvalue of $M$ will then be controlled by the product of the largest matrix entry of $\boxed{M}$, and the entrywise $L^1$-norm of the family of matrices $M'$. 

For the $L^1$-norm, we inspect the entries in the matrix $M'$. We know $r\approx \eta$, and since $\v_\epsilon$ $C^0$-approximates $\v$ we obtain $\v_\epsilon \approx 1$. So to get a uniform upper bound on the $L^1$-norm of $M'$, we just need to make sure our smoothing is chosen so that:
\begin{itemize}
    \item $\v_\epsilon '$ is bounded above by a uniform constant.
\end{itemize}

Concerning the matrix entries of $\boxed{M}$, we would like to pinch those close to zero. Since $r\approx \eta$ and $\v_\epsilon \approx 1$, this will require:
\begin{itemize}
    \item $\v_\epsilon '$ can be pinched arbitrarily close to zero,
    \item $\v_\epsilon ''$ can be pinched arbitrarily close to zero.
\end{itemize}

Consolidating all the required conditions from the bulleted points above, we see that it suffices to choose our smoothings $\v_\epsilon$ so that both $|\v_\epsilon'|$ and $|\v_\epsilon ''|$ are arbitrarily small, e.g. both $< \e/66$ suffices. The desired smoothing is provided by Lemma \ref{lem:smoothing at r=eta} in the Appendix.

\medskip

\medskip

\noindent{\underline{Singularity at $r=\delta$}:} In the vicinity of $r=\delta$, the functions $f$ and $\v$ are given by:
\begin{equation*}
f(r)=\begin{cases}
r &\d<r\\
\d &r\le \d
\end{cases},
\qquad\qquad
\v(r)=\begin{cases}
h(r)&\d <r \\
0 & 0<r \le \d
\end{cases}
\end{equation*}
Since both functions are $C^0$, but fail to be $C^1$ at $r=\delta$, they will both need to be smoothed. Also, observe that the smoothing $f_\epsilon$ will have $f_\epsilon '$ that transitions from $0$ to $1$. While the smoothing $\v_\epsilon$ has $\v_\epsilon '$ transitioning from $0$ to $h'(\delta)\approx \frac{1}{\delta \ln{(\eta/\delta)}}$. Note that as $\delta \to 0^+$, we have $h'(\delta) \to +\infty$, so decreasing $\delta$ does not make $\v$ any ``closer'' to a $C^1$-function. 

We can again analyze the coordinate plane curvature Equations (\ref{eqn:curvature component})-(\ref{eqn:curvature component 2}), which we want to be $\leq \e/22$. Arguing as in the previous cases, this reduces to choosing the smoothings $f_\epsilon, \v_\epsilon$ so that the following constraints hold:
\begin{itemize}
    \item $\v_\epsilon' f_\epsilon' \v_\epsilon f_\epsilon \geq -\e/44$, used for Equations (\ref{eqn:curvature component 4}) and (\ref{eqn:curvature component 5}),
    \item $\v_\epsilon^2 \cdot \left((f_\epsilon ')^2-1\right)\geq -\e/44$, used for Equations (\ref{eqn:curvature component 4}) and (\ref{eqn:curvature component 5}),
    \item $\v_\epsilon f_\epsilon\left( \v_\epsilon ''f_\epsilon+2\v_\epsilon'f_\epsilon'+\v_\epsilon f_\epsilon''\right)\geq -\e/22$, used for Equations (\ref{eqn:curvature component 6}) and (\ref{eqn:curvature component 7}),
    \item $\v_\epsilon ^2f_\epsilon^2\left(\v_\epsilon^2-(\v_\epsilon 'f_\epsilon+\v_\epsilon f_\epsilon')^2\right)\leq \e/22$, used for Equation (\ref{eqn:curvature component 8}),
    \item $\fe'' \geq 0$, used for Equation \eqref{eqn:curvature component 2}.  
\end{itemize}
Concerning the control of the eigenvalues of the matrix $M:=M'+\boxed{M}$, this breaks down into bounding the $L^1$-norm of $M'$ and pinching the entries of $\boxed{M}$. Looking at the entries of $M'$, and recalling that $|\v_\epsilon|$ is uniformly bounded, the $L^1$-norm condition translates to
\begin{itemize}
    \item $|\v_\epsilon'f_\epsilon + \v_\epsilon f_\epsilon'|$ is uniformly bounded above
\end{itemize}
while looking at the entries of $\boxed{M}$, pinching the entries reduces to choosing a smoothing so that the following expressions are arbitrarily close to zero
\begin{itemize}
    \item $|\v_\epsilon' f_\epsilon' \v_\epsilon f_\epsilon|$, 
    \item $\v_\epsilon^2 \cdot |(f_\epsilon ')^2-1|$, 
    \item $\big|\v_\epsilon f_\epsilon\left( \v_\epsilon ''f_\epsilon+2\v_\epsilon'f_\epsilon'+\v_\epsilon f_\epsilon''\right)\big|$, and
    \item $|\v_\epsilon ' \v_\epsilon f_\epsilon|$.
\end{itemize}

Now recall that we can choose our approximating functions arbitrarily close in the $C^0$-metric, which implies the approximations $f_\epsilon\approx \delta, \v_\epsilon \approx 0$ can be chosen as close as we want. In view of this, and consolidating the bulleted constraints above, we see that it suffices to ensure that our smoothings $f_\epsilon, \v_\epsilon$ satisfy the following conditions:
\begin{enumerate}
    \item $|f_\epsilon'|$ is uniformly bounded, $f''_\epsilon\ge 0$ and $f_\epsilon''=O\left(\frac{1}{\d}\right)$;
    \item $f_\epsilon \v _\epsilon' = O\left(\frac{1}{\ln{(\eta/\d)}}\right)$;
    \item $f_\epsilon^2 \v _\epsilon '' = O\left(\frac{1}{\ln{(\eta/\d)}}\right)$.
\end{enumerate}
Note that the last two conditions are {\it mixed} conditions, in that they require some interdependence between the two smoothings $f_\epsilon, \v_\epsilon$. In particular, the two smoothings cannot be chosen independently. In conditions (1)-(3), the implicit constants have to be universal. 

\medskip

Finally, we note that for each of the three singularities, we have reduced control of the curvatures to establishing smoothings of the functions $f$, $\phi$ with certain prescribed 
properties. In the Appendix, we construct the smoothings, see Lemmas \ref{lem:smoothing at r=0}, \ref{lem:smoothing at r=eta}, and \ref{lem:smoothing at r=delta}. These smoothings depend on the parameter $\delta>0$, and choosing $\delta$ close enough to $0$ ensures the smoothings have the desired control. This 
establishes {\bf Claim 3}.

\medskip

\medskip

To complete the proof of the theorem, we are left with controlling the volume of $(N, g)$. 
We will denote by $\Vol(S)$ the volume of $S$ with respect to the locally symmetric metric $\h_{m-1} + \h_{n-1}$. 
We first estimate the volume for the singular metric in each of the regions $r\leq 0$, $0\leq r\leq \d$, $\d \leq r \leq \eta$, near $S$. Note that in the range $r\geq \eta$ the singular metric (and hence its volume) coincides with the locally symmetric one.
In the range $r\leq 0$, the singular metric reduces to the local product metric $\h_{m-1} + \h_{n-1} + (\d e^{\d r})^2 d\theta^2+dr^2$. Thus on that portion of the space the volume is exactly
$$2\pi \cdot\Vol (S)\cdot \int_{- \infty}^0 \d e^{\d r} \, dr =2\pi \cdot\Vol (S).$$ When $0\leq r \leq \d$, the singular metric is again a local product metric $\h_{m-1} + \h_{n-1} + \d^2 d\theta^2+dr^2$, which gives a cylinder of radial distance $r$ about $S$, hence has volume $2\pi \d^2 \cdot \Vol(S)$. Finally, within the region $\d\leq r\leq \eta$, the identity map from the locally symmetric metric to our warped metric is $1$-Lipschitz (recall that $h(r)\leq 1$). So on this region, the volume is bounded above by the locally symmetric one. 
Summing up these inequalities about each region near $S$ establishes that, for the singular metric we have $\Vol(N) < V$ 
for some uniform constant $V$. From our estimates, the constant $V$ can be chosen to be any real number that is greater than 
$$\Vol(M)+ 2\pi \cdot \Vol(S).$$ Finally, the smoothing only pointwise perturbs the metric by amounts that are $O(\delta)$, so when $\delta$ is chosen small enough, the smoothed out metric still satisfies the same estimate $\Vol(N,g) < V$.
This completes the proof of Theorem \ref{thm:main theorem 2} (1).

\end{proof}

\subsection{Almost non-positively curved metrics on branched covers}
To prove Theorem \ref{thm:main theorem 2} (2) we will need the following result from \cite{fornari1990ramified}.

\begin{theorem}[\cite{fornari1990ramified} Section 2]\label{thm:FS}
Let $(M,g_M)$ be a non-positively curved manifold, $(S,g_S)$ an embedded totally geodesic, codimension two submanifold of $M$, and suppose that, on some $r$-tube about $S$, we have that $g_M = g_S + e_2$.  
Let $X_d$ denote the $d$-fold cyclic branched cover of $M$ about $S$, and $\phi_d : X_d \to M$ the associated ramified covering map.  
Then, given any $\d_1, \d_2 > 0$ with $0 < \d_1 < \d_2 < r$, there exists a smooth, non-positively curved metric $g$ on $X_d$ such that
\begin{enumerate}
\item $g = g_M = g_S + e_2$ on the $\d_1$-neighborhood of $S$, and
\item $g = \phi_d^* g_M$ outside of the $\d_2$-neighborhood of $S$.  In particular, on the $(\d_2, r)$-annulus of $S$, we have $g = g_S + \phi_d^*(e_2)$.  
\end{enumerate}
\end{theorem}

Note that, in Theorem \ref{thm:FS}, we abuse notation and also use $S$ to refer to the ramification locus of $X_d$.  

\begin{remark}\label{rmk:necessity of product lattice}
(1) The smoothing argument carried out in \cite{fornari1990ramified} is done locally near the branching locus $S \subset X_d$. They show how to modify the pull-back metric $\phi_d^* g_M$ within the $\d_2$-neighborhood of the branching locus $S$ to ``flatten out'' the singularity. Moreover, it can be carried out so that the volume of the smoothed out Riemannian metric is as close as desired to the volume of the singular pull back metric.

\noindent (2) The use of Theorem \ref{thm:FS} is essential in Fornari and Schroeder's proof of Theorem \ref{thm:main theorem 1} (2) in \cite{fornari1990ramified}. The natural pull back of the Riemannian metric to the branched cover yields a metric with singularities transverse to $S$, due to the excess angle arising from the branched covering construction. One can then attempt to smooth out this singularity. Theorem \ref{thm:FS} shows this can be done when the metric splits off a transverse flat factor. This allows a different approach to the smoothing problem, as one can instead try to construct a metric on the original $M$ which satisfies the local transverse flatness condition near $S$, and then appeal to Theorem \ref{thm:FS} to smooth the lifted metric on $X_d$.

\noindent (3) In the case where $(M,S)$ is modeled on $(\H^m \times \H^n, \H^{m-1} \times \H^{n-1})$, one can edit the argument of Fornari and Schroeder so that all curvatures of this metric lie within $(-1-\e, 0]$ for any prescribed $\e > 0$.  
Given some large positive constant $k_0$, their argument\footnote{In \cite{fornari1990ramified}, the authors use $\tilde{\rho}$ for $r$ and $k$ for the branching degree $d$.} relies on constructing a smooth function $\v(r)$ that varies from $f(r) = r$ to $h(r) = d(r-k_0)$.
One can do this in a convex manner so that all curvatures are non-positive.  
If one just smooths out the intersection point of these functions, then $\v(r)$ will have a large second derivative near this point and thus the metric will have large negative curvature.  
But, for $k_0$ chosen sufficiently large, one can construct $\v(r)$ in such a way that $\v'' < \e$ for all $r$.  
This keeps all curvatures bounded below by $-1-\e$.  
This point is only needed for Remark \ref{rmk:pinched curvature} below.
\end{remark}

We can now use Theorem \ref{thm:FS} and an essentially identical argument as in the proof of Theorem \ref{thm:main theorem 2} (1) to prove Theorem \ref{thm:main theorem 2} (2).  

\begin{proof}[Proof of Theorem \ref{thm:main theorem 2} (2)]
In view of Theorem \ref{thm:FS}, it suffices to construct a metric $g$ on $M$ that is (1) non-positively curved within a small tubular neighborhood of $S$, (2) is almost non-positively curved away from this small tubular neighborhood, and (3) satisfies the ``local transverse flatness'' condition within the tubular neighborhood of $S$. Theorem \ref{thm:FS} then provides a smoothing of the lifted metric on $X_d$ in the small tubular neighborhood of the branching locus $S \subset X_d$ while maintaining non-positive curvature on the tubular neighborhood. Thus the resulting smoothed metric on $X_d$ will have the same small positive upper bound on sectional curvatures as the original metric on $M$.

One can essentially construct such a metric on $M$ by just removing the cusp  portion (corresponding to $r < 0$) from the metric constructed in the proof of Theorem \ref{thm:main theorem 2} (1).
In particular, for $\d, \eta > 0$ chosen as in the proof of Theorem \ref{thm:main theorem 2} (1), we define $\v \colon [0,\infty) \to \R$ by
\begin{equation*}
\v(r)=\begin{cases}
1 & r \ge \eta\\
\frac{1}{\ln \left( \frac{\eta}{\d} \right) } \ln \left( \frac{r}{\d} \right) &\d <r < \eta\\
0 & 0 \leq r \le \d 
\end{cases},
\end{equation*}

We define the metric $g'$ via Equation \eqref{eqn:polar coordinate warped product metric} by defining $f$, $\a$, and $\b$ as
\begin{equation*}
f(r) = r, \qquad \a(r, \th) = \cosh(\v(r) r \cos(\th)), \qquad \b(r, \th) = \cosh(\v(r) r \sin(\th)).
\end{equation*}
Note that, for $0 \leq r < \d$, we have
\begin{equation*}
\a(r, \th) = \cosh(0 \cdot r \cos(\th)) = 1, \qquad \qquad \b(r, \th) = \cosh(0 \cdot r \sin(\th)) = 1.
\end{equation*}
Thus, for $0 \leq r < \d$, from \eqref{eqn:polar coordinate warped product metric} we have
\begin{equation*}
g' = \h_{m-1} + \h_{n-1} + r^2 d \th^2 + dr^2 = \h_{m-1} + \h_{n-1} + e_2.
\end{equation*}

The metric $g'$ defined here is not equivalent to the metric defined in the proof of Theorem \ref{thm:main theorem 2} (1).  
On the interval $(0,\d)$ we have here $f(r) = r$, whereas in the previous argument we had $f(r) = \d$.  
But one easily checks that all of the estimates in the proof of Theorem \ref{thm:main theorem 2} (1) still apply since $\v \equiv 0$ on this region.  
Thus, one may choose $\d, \eta > 0$ sufficiently small so that the sectional curvature of $g'$ is less than $\e/2$ where $r\neq \d$ and $r\neq \eta$.  
One can then smooth $g'$ to obtain a Riemannian metric $g$ defined on $M$ keeping all sectional curvatures less than $\e$. Near $r=\eta$, the smoothing $\varphi_\epsilon$ is exactly the same as in the proof of Theorem \ref{thm:main theorem 2}(1) and the same calculations yield the same curvature bound. Now near $r=\d$, we use the $\varphi_\epsilon$ constructed in Lemma \ref{lem:smoothing at r=delta}, while using $f_\epsilon=f=r$. To see the constraints in the proof of Theorem \ref{thm:main theorem 2}(1) still hold with this choice of $f_\epsilon$, we need to ensure that the following conditions are satisfied:
\begin{itemize}
    \item $|f_\epsilon'|$ is uniformly bounded, $f_\e''\geq 0$, and $f_\epsilon'' = O(\frac{1}{\delta});$
    \item $r \v _\epsilon' = O\left(\frac{1}{\ln{(\eta/\d)}}\right)$;
    \item $r^2 \v _\epsilon '' = O\left(\frac{1}{\ln{(\eta/\d)}}\right)$.
\end{itemize}

Since $f_\epsilon(r)=r$, the first property is immediate. The last two properties follow easily from the construction of $\v_\epsilon$ in Lemma \ref{lem:smoothing at r=delta}, along with the fact that $r<\d+\epsilon$ (these are the same estimates that establish properties (7) and (8) in the proof of Lemma  \ref{lem:smoothing at r=delta}).

Observe that our smoothed out metric $g$ will satisfy all three conditions (1),(2), and (3) that we stated in the beginning of this proof. Applying Theorem \ref{thm:FS}, we obtain a smooth metric $g_d$ on $X_d$ with all sectional curvatures bounded above by $\e$.  

The last thing that we need to consider is the volume of $(X_d, g_d)$.  
Note that, when defining $g'$ above,
\begin{equation*}
\a(r, \th) \leq \cosh(r \cos(\th)) \qquad \b(r, \th) \leq \cosh(r \sin(\th))
\end{equation*}
for all values of $r$ and $\th$. Moreover, after smoothing, $f=r$ and $0\leq\v_\epsilon\leq 1$. Therefore our smooth metric $g$ satisfies $\text{Vol}(M,g) \leq \text{Vol}(M, \h_m + \h_n)$. This tells us the pull-back singular metric $\phi_d^*g$ on $X_d$ also satisfies the inequality 
\begin{equation*}
\text{Vol}(X_d, \phi_d^*g) \leq d \cdot \text{Vol}(M, \h_m + \h_n)
\end{equation*}
From Remark \ref{rmk:necessity of product lattice} (1), the smoothing on $X_d$ can now be carried out to have volume as close to that of $\text{Vol}(X_d, \phi_d^*g)$ as desired. And hence, can be bounded above by $d \cdot \text{Vol}(M, \h_m + \h_n) +\xi$ for any prescribed $\xi>0$.
\end{proof}

\begin{remark}\label{rmk:pinched curvature}
For the same reasons as in the proof of Theorem \ref{thm:main theorem 2} (1), the metric constructed in the above proof will contain very large negative sectional curvatures.  
The lower sectional curvature bound will approach $- \infty$ as $\e \to 0$.  
But, if the normal injectivity radius of $S$ in $M$ is sufficiently large, then one may construct a metric on $X_d$ whose sectional curvatures lie in the interval $(-1-\e , \e)$ and with the same volume bound.  
For large positive constants $A$ and $B$, one defines $\v$ to be a smooth, increasing function which is identically equal to $0$ for $r \leq A$, and is equal to $1$ for $r \geq B$.  
One then defines $g'$ exactly as in the above argument.  
One can check that, for $A$ chosen sufficiently large, all of the estimates in the proof of Theorem \ref{thm:main theorem 2} (1) can be forced arbitrarily close to $0$.  
One then obtains the desired metric on $X_d$ by applying Remark \ref{rmk:necessity of product lattice} (3).  
\end{remark}

\section{The special case of product lattices}\label{sect:product case}

In this section, we focus on the particular case of product lattices. This means that $M$ splits as a product, and the connected submanifold $S$ also splits as a product. In this special case, we can establish stronger statements, both geometrically and algebraically.

\subsection{Geometric properties} In this subsection we establish the proof of Theorem \ref{thm:main theorem 1}. We start by recalling some definitions.

\begin{definition}\label{def:A-regular}
    A Riemannian metric $g$ with curvature tensor $R$ is {\it $A$-regular} if there exists a sequence of positive numbers $\{ A_k \}$ such that, for each $k \geq 0$, the $k^{th}$ covariant derivative of $R$ satisfies $| \nabla^k R|_{C^0} < A_k$.  
\end{definition}

    Note that when $k = 0$, this implies that the sectional curvature of $g$ is bounded.  
    Obviously, any metric on a compact manifold is $A$-regular, so this notion is only meaningful for non-compact manifolds. It is easy to see that symmetric spaces are $A$-regular. Therefore, in our setting, we just need to show that our constructed metric is $A$-regular on the region about $S$ where we turn $S$ into a cusp.

\begin{proof}[Proof of Theorem \ref{thm:main theorem 1} (1)]
Let $M$ be a product $M_1 \times M_2$, and $S = S_1 \times S_2$ with $S_i \subset M_i$ a connected embedded totally geodesic submanifold. We want to put a complete finite-volume $A$-regular non-positively curved metric on $N= M\setminus S$. As before, we do this by locally modifying the metric to produce cusps near $S$.

As a first step, we apply Proposition \ref{prop:product metrics} to obtain a metric $g_2$ on $M$ which is locally flat transverse to $S$. Let $\eta$ denote the normal injectivity radius of $S$, and let $\eta^\prime<\eta$ chosen small enough so that in the $\eta^\prime$-neighborhood of $S$, the metric $g_2$ is locally isometric to $\h_{m-1} \oplus \h_{n-1} \oplus e_2$. Then the universal cover of that neighborhood can be identified with the $\eta^\prime$-neighborhood of $\H^{m-1} \times \H^{n-1}\times \{0\}$ inside $\H^{m-1} \times \H^{n-1}\times \R^2$. As usual, we will identify hyperbolic space $\H^k$ with $\R^k$ equipped with the hyperbolic metric $\h_k$. 

We now define a metric $\tilde{g}_1$ on $\R^{m-1}\times\R^{n-1}\times\S^1\times(-\infty,\e_1)$ where $\e_1=\frac{1}{2}\eta'$. 
From Equation \eqref{eqn:polar coordinate warped product metric}, by setting $\alpha=1$ and $\beta=1$, we have a metric 
\begin{align}
    \tilde{\l}=\tilde{\h}_{m-1}+\tilde{\h}_{n-1}+f^2(r)d\th^2+dr^2.
\end{align}

We want to construct a function $f(r)$ such that $f'> 0$, $f''\ge 0$, and $f \to 0$ exponentially as $r \to - \infty$. 
Let $\e_2=\frac{1}{2}\e_1$.
Then there exists $c<\text{min}\{1,\e_2\}$ and an $\epsilon>0$ so that $\epsilon e^c=c$. 
Let $\frac{1}{2}c=\e_3$ and define two functions $f_1(r)$ and $f_2(r)$ by
\begin{align}
    f_1(r)=\epsilon e^r,\quad & r\in[\e_3,c],\\
    f_2(r)=r,\quad & r\in[c,\e_2].
\end{align}

Notice that $f_1(c)=\epsilon e^c=c=f_2(c)$ and $f'_1(c)=\epsilon e^c=c\leq 1=f'_2(c)$. 
Therefore, by Lemma \ref{lem:smoothing lemma}, there exists a smooth function $f$ on $[\e_3,\e_2]$ such that $f'' \geq 0$,  $f_1(\e_3)=f(\e_3)$, and $f_2(\e_2)=f(\e_2)$. 

We thus have a smooth metric $\tilde{g}_1$ defined on $\R^{m-1}\times\R^{n-1}\times\S^1\times(-\infty,\e_1)$ given by
\begin{align}\label{eqn:flat-to-cusp}
    \tilde{g}_1&=\begin{cases}
        \tilde{\h}_{m-1}+\tilde{\h}_{n-1}+(\epsilon e^r)^2d\th^2+dr^2, \quad & r\in (-\infty,\e_3]\\
        \tilde{\h}_{m-1} + \tilde{\h}_{n-1} + f^2(r) d \th^2 + dr^2, \quad &r \in [\e_3, \e_2] \\
        \tilde{\h}_{m-1}+\tilde{\h}_{n-1}+r^2d\th^2+dr^2, \quad & r\in [\e_2,\e_1]
    \end{cases}.
\end{align}

By Corollary \ref{cor: only nonzero term}, the sectional curvature of $\tilde{g}_1$ is non-positive. 
From our construction, it is immediate that the metric $\tilde{g}_1$ is also invariant under isometries of $\H^{m-1}\times \H^{n-1}$, and under rotations and reflections in the $\R^2$ factor. Thus it descends to a well-defined metric $g_1$ on the (punctured) tubular neighborhood of $S$. Note that, on the annular region $\R^{m+n-2} \times \S^1 \times [\e_2, \e_1]$, the metric $\tilde{g}_1 =  \tilde{g}_2$, as they both coincide with the metric $\tilde{\h}_{m-1} + \tilde{\h}_{n-1} + e_2$. This allows us to glue the $g_1$ metric on the inside of the punctured tubular $r\in(-\infty, \e_1)$ region, with the $g_2$ metric on the outside of the smaller $\e_2$-neighborhood of $S$.
These metrics glue together to yield a non-positively curved metric $g$ defined on $M\setminus S$.

Since $\a = \b = 1$ and $f = \epsilon e^r$ on $(- \infty, 0)$, by Remark \ref{rmk:finite-volume} the metric $g$ has finite-volume.
For $A$-regularity, consider the metric $g$ near the cusp. 
By Corollary \ref{cor: only nonzero term} there are only three nonzero terms of the curvature tensor, all identically equal to $-1$.  
Therefore all of its derivatives are zero and hence bounded. 
Thus the metric $g$ is $A$-regular, as desired. 
\end{proof}

\begin{remark}
    It is possible to extend the proof to some cases where $S$ is not connected. As long as the projection of $S$ to the two factors consists of a collection of pairwise disjoint codimension one submanifolds. Then $S \subset Y_1\times Y_2$, where $Y_i \subset X_i$. One can then perform this same modification locally in a tubular neighborhood of each of the finitely many connected components of $Y_i$, just choosing the $\eta_3$ small enough to guarantee the tubular neighborhoods are pairwise disjoint. The resulting product metric on $M$ will be locally flat transverse to each connected component of $S$, and the rest of the argument goes through unchanged.
\end{remark}

\subsection{Group theoretic properties}\label{sect:proof of group properties}

In this section, we will establish Theorem \ref{thm: group theoretic properties}, which concerns various group theoretic properties of $\pi_1(N)$ in the particular case where $M$ is a compact product, and $S$ is connected. In this case, the group $\pi_1(N)$ has some graph of groups decompositions, which one can exploit. In the irreducible lattice case, and in the case where $S$ is not connected, this tool is no longer available.

Let us start by describing the graph of groups splittings that we will use.  In the product lattice case, recall that $N=(X_1\times X_2) \setminus S$. Recall that $S$ projects to a pair of codimension one totally geodesic submanifolds $Y_i \subset X_i$, and $S=Y_1 \times Y_2$. There are now three slightly distinct possibilities that arise, according to whether the $Y_i$ separate $X_i$. Denote by $W_i$ the open manifold $X_i \setminus Y_i$. We will call $N$ 
\begin{itemize}
    \item Type I, if $W_1$ and $W_2$ are both connected;
    \item Type II, if $W_1$ is connected, and $W_2$ is disconnected;
    \item Type III, if $W_1$ and $W_2$ are both disconnected.
\end{itemize}
Up to relabeling the two factors, these exhaust our possibilities. The graph of group descriptions one obtains varies slightly in these three cases. 

\begin{lemma}\label{lem:graph-of-group-codim-zero}
    In the product case where $M$ is compact and $S$ is connected, the group $\pi_1(N)$ has a finite graph of group decomposition with the property that
    \begin{itemize}
        \item vertex groups are fundamental groups of a product of a closed hyperbolic manifold and a compact hyperbolic manifold with totally geodesic boundary;
        \item edge groups are fundamental groups of a product of two compact hyperbolic manifolds with boundary.
    \end{itemize}
    Moreover each edge group sits as a quasi-convex subgroup in the adjacent vertex groups.
\end{lemma}

\begin{proof}
    The graph of group decompositions are obtained by applying the Seifert-van Kampen theorem to $N$. Indeed, the manifold $N=(X_1\times X_2) \setminus S$ is naturally covered by the two open sets $W_1\times X_2$ and $X_1\times W_2$ (recall $W_i:=X_i\setminus Y_i$). These open sets each have one or more connected components, according to the type of $N$. Moreover, note that the (intrinsic) metric completion $\overline{W_i}$ of $W_i$ is a compact hyperbolic manifold with totally geodesic boundary, and $\pi_1(W_i) \cong \pi_1(\overline{W_i})$. The boundary components of $\overline{W_i}$ are each isometric to $Y_i$. Applying Seifert-van Kampen to this covering, we see that the vertex groups are of the form $\pi_1(\overline{W_1}^0 \times X_2)$ or $\pi_1(X_1 \times \overline{W_2}^0 )$, where $\overline{W_i}^0$ denotes a connected component of $\overline{W_i}$.

    For the edge groups, we need to look at the connected components of the intersection of the two sets. These groups are of the form $\pi_1(\overline{W_1}^0 \times \overline{W_2}^0)$, as desired. Finally, the maps from the edge groups to the vertex groups are induced by an inclusion of spaces. Since the subspace corresponding to the edge groups are convex subsets in the spaces corresponding to the vertex group, applying the $\pi_1$-functor gives rise to quasi-convex subgroups.
\end{proof}

With Lemma \ref{lem:graph-of-group-codim-zero} in place, we can now establish Theorem \ref{thm: group theoretic properties}, which we restate for convenience.

\begin{theorem}
    
Assume $M = X_1\times X_2$ is a compact product, and $S\subset M$ a connected totally geodesic codimension two modeled on $\H^{m-1}\times \H^{n-1}$. Then for
the complement $N=M\setminus S$, the fundamental group $\pi_1(N)$ satisfies the following properties:
\begin{enumerate}
\item uniform exponential growth,
\item Tits alternative,
\item Kazhdan subgroups are always trivial,
\item finite asymptotic dimension,
\item Baum--Connes conjecture with coefficients,
\item property RD (rapid decay).
\end{enumerate}
\end{theorem}  

\begin{proof}
All of the properties in our theorem can be deduced from the graph of groups decompositions of $\pi_1(N)$ described in Lemma \ref{lem:graph-of-group-codim-zero}.  

\vskip 5pt

\noindent \underline{1. Uniform Exponential Growth.} The edge groups in the graph of groups decomposition in Lemma \ref{lem:graph-of-group-codim-zero} have infinite index in the vertex groups, which immediately implies uniform exponential growth, see Bucher--de la Harpe \cite{BdlH}.

\vskip 5pt

\noindent \underline{2. Tits Alternative.} We use the graph of groups decomposition for $\pi_1(N)$ given in Lemma \ref{lem:graph-of-group-codim-zero}. The vertex groups are the fundamental group of a compact non-positively curved manifold with totally geodesic boundary. Each vertex group is a direct product of a lattice in hyperbolic space (from the $X_i$-factor), with the fundamental group of a hyperbolic manifold with totally geodesic boundary (from the closure of the connected component in the $W_i$-factor). As such, each vertex group is a retract of a product uniform lattice in the product Lie group $\SO(m,1) \times \SO(n,1)$, so satisfies the classical Tits alternative. Thus any subgroup of a vertex group either contains a non-abelian free group, or is virtually solvable. We can in fact recover a bit more information about the solvable subgroups. If $H$ is such a solvable subgroup, we can project it onto the two factors obtaining subgroups $H_1$, $H_2$ inside a uniform torsion-free lattice subgroup in $\SO(m,1)$, $\SO(n,1)$. These lattices are Gromov hyperbolic and torsion-free, so their solvable subgroups must be infinite cyclic. This tells us $H$ is contained in a $\mathbb Z^2$-subgroup, hence is itself abelian (of rank at most two). We have shown that the vertex groups in the graph of groups decomposition satisfy a strong form of the Tits alternative: subgroups are either virtually abelian of rank at most two, or contain a non-abelian free subgroup. Applying \cite[Prop 6.18]{FLS} allows us to conclude that $\pi_1(N)$ also satisfies the Tits alternative. In fact, the proposition also gives us the slightly stronger conclusion that any solvable subgroup in $\pi_1(N)$ has derived length $\leq 3$.

\vskip 5pt

\noindent \underline{3. Kazhdan subgroups.} Let $H$ be a Kazhdan subgroup of $\pi_1(N)$. Kazhdan groups have Property (FA). Therefore $H$ has to fix some vertex in the Bass--Serre tree associated to the splitting of
$\pi_1(N)$, and hence is conjugate to a subgroup of a vertex stabilizer. So we can view $H$ as isomorphic to a Kazhdan subgroup of one of the vertex groups.
From the description of vertex groups in Lemma \ref{lem:graph-of-group-codim-zero}, the vertex group splits naturally as a direct product, and inherits a product action on $\mathbb H^m\times \mathbb H^n$.
The first factor acts fixed point freely and isometrically on $\mathbb H^m$, so from property (FH) we have that $H$ must lie in the kernel of the 
projection map, see Bekka--de la Harpe--Valette \cite[Prop 2.6.5]{BdlHV}). Applying the same reasoning to the projection onto the second factor, we conclude $H$ is trivial.

\vskip 5pt

\noindent \underline{4. Finite asymptotic dimension.} Bell and Dranishnikov \cite{Bell-Dranishnikov} have shown that the fundamental group of a graph of groups has finite asymptotic dimension provided all the vertex groups have finite asymptotic dimension. From the description of the vertex groups in Lemma \ref{lem:graph-of-group-codim-zero}, we know these are products of Gromov hyperbolic groups. It is well-known that Gromov hyperbolic groups have finite asymptotic dimension (see e.g. Roe \cite{Roe}), and that products of groups with finite asymptotic dimension still have finite asymptotic dimension (see e.g. \cite[Thm 32]{Bell-Dranishnikov-survey}). Thus vertex groups have finite asymptotic dimension, and applying Bell--Dranishnikov property (4) follows.

\vskip 5pt

\noindent \underline{5. Baum-Connes Conjecture.} It is a result of Oyono-Oyono \cite{Oyono} that the Baum-Connes Conjecture (with coefficients) holds for a graph of groups if and only if it holds for the vertex groups. In our graph of groups decomposition for $\pi_1(N)$, the vertex groups are direct products of a uniform lattice in $\SO(m,1)$ (or $\SO(n,1)$), with a retract of a lattice in $\SO(m,1)$ (or $\SO(n,1)$). The Baum-Connes Conjecture holds for lattices and retracts of lattices in the Lie groups $\SO(m,1)$, $\SO(n,1)$, see e.g. Lafforgue \cite{Lafforgue2}, and is preserved under direct products of groups, by Chabert--Echterhoff \cite{Chabert}. Applying Oyono-Oyono, we obtain the desired result.

\vskip 5pt

\noindent \underline{6. Property (RD).}
Inspired by Haagerup  \cite{Haagerup}, Jolissaint introduced and studied the Rapid Decay Property (RD) for groups \cite{Jolissaint} (see \cite{Chatterji-intro} for a gentle introduction to this topic). It has since been used to prove certain cases of the Novikov and Baum-Connes Conjecture \cite{Lafforgue}. Recent work by Chatterji--Gautero shows that, if a graph of groups has vertex groups satisfying property (RD), and if the edge groups have loose polynomial distortion in the vertex groups, then the entire group satisfies property (RD) (see \cite[Prop 1.3]{Chatterji-Gautero}). For the vertex groups, we have from Lemma \ref{lem:graph-of-group-codim-zero} that they are products of Gromov hyperbolic groups. But property (RD) is known to hold for Gromov hyperbolic groups (see \cite{dlH}), and to be inherited under direct products (see \cite{Jolissaint}). Moreover, arguments similar to those in Lemma \ref{lem:graph-of-group-codim-zero} easily show the edge groups are indeed linearly (and hence polynomially) distorted in the vertex groups. Applying Chatterji--Gautero completes the proof of property (RD).
\end{proof}

\begin{remark}
    In Theorem \ref{thm: group theoretic properties}, the hypothesis that $S$ is connected is used in a couple of ways. First, it guarantees that the projections $Y_i$ of $S$ onto the two factors are {\it embedded} codimension one submanifolds. If $S$ has multiple connected components, this is no longer guaranteed. Secondly, even if the $Y_i$ are embedded, it could be possible in the disconnected case that $S \neq Y_1\times Y_2$. In that situation the sets described in Lemma \ref{lem:graph-of-group-codim-zero} would not cover $N$, and hence would not give rise to a graph of groups structure for $\pi_1(N)$. If there are multiple connected components for $S$, as long as $S = Y_1\times Y_2$, the arguments in Theorem \ref{thm: group theoretic properties} still go through.
\end{remark}

\section*{Appendix: Smoothing continuous functions with controlled derivatives}\label{Section:Introduction}

In this appendix we collect some results about smoothing continuous functions that were used previously in the paper.  
The first proposition below is basically a restatement of \cite[Proposition 13.5]{belegradek2012complex}.  

\begin{proposition}\label{lem:smoothing lemma}
Suppose that $f \colon [a,b] \to \R$ is a positive continuous function such that, for some $c \in (a,b)$, the restrictions of $f$ to $[a,c]$ and $[c,b]$ are smooth with $f'' \geq k$ and $f'(c_{-}) \leq f'(c_{+})$.  
Then for any $\epsilon > 0$ sufficiently small, there exists a smooth function $f_\epsilon \colon [a,b] \to \R$ which satisfies
\begin{enumerate}
    \item $f_\epsilon = f$ outside of the $\epsilon$-neighborhood of $c$.
    \item $f_{\epsilon}'' \geq k$.
    \item if $f$ is non-decreasing, then $f_{\epsilon}' \geq 0$.
\end{enumerate}
\end{proposition}

Proposition 13.5 of \cite{belegradek2012complex} states this result for $f$ increasing and $f_\epsilon' > 0$, but the same argument holds for non-decreasing and $f_\epsilon' \geq 0$.  
The function $f_\d$ is constructed using convolution as described, for example, in \cite[Chapter 2]{Hirsch}.
If $f$ is non-decreasing, then convolution with a non-decreasing {\it convolution kernel} produces a non-decreasing function.  
Thus, $f_{\epsilon}' \geq 0$.  

The remaining three lemmas are all used to smooth functions in the proof of Theorem \ref{thm:main theorem 1} (1).  
We present the lemmas in the order that they appear in the argument.
AI was used to construct explicit functions which satisfy the properties of each lemma, and the authors verified the correctness of these solutions.  
The first two lemmas, which smooth one function at either $r=0$ or $r=\eta$, are reasonably clear from just drawing a picture.  
So we only present the construction from AI for the smoothing at $r = \d$, which is considerably more complicated due to the required relationships between $\fe$ and $\ve$.  

In what follows, $\d , \eta > 0$ are constants with $\d \ll \eta$.  

\begin{lemma}[Smoothing at $r=0$]\label{lem:smoothing at r=0} 
Let $\e > 0$ be arbitrary.  
Then there exists $\d > 0$ sufficiently small and a smooth function $\fd \colon \R \to \R$ which satisfies the following properties:
\begin{enumerate}
\item $\fd(r) = \d e^{\d r}$ for all $r \leq -1$.  
\item $\fd(r) = \d$ for all $r \geq 0$.  
\item $|\fd(r) - \d e^{\d r}| < \e$ for all $r \in (-1, 0)$.  
\item $- \frac{\fd''(r)}{\fd(r)} < \e $ for all $r \in (-1, 0)$.

\end{enumerate}
\end{lemma}

\begin{proof}[Proof of Lemma \ref{lem:smoothing at r=0}]
Define a function $\fd(r)$ by $$\fd(r)=\d(1+\chi(r)(e^{\d r}-1)),$$
where $\chi(r):\mathbb{R}\to[0,1]$ is a smooth non-increasing function such that $$\chi(r)=1 \text{ for } r\leq -1, \text{ and } \chi(r)=0 \text{ for } r\geq 0.$$
Then this function directly satisfies the conditions $(1)-(3)$. Since $|\frac{\fd''}{\fd}|=O(\d)$ on $(-1,0)$, the property $(4)$ holds for sufficiently small $\d>0$.
\end{proof}

\begin{lemma}[Smoothing at $r=\eta$]\label{lem:smoothing at r=eta} 
Let $\e > 0$ arbitrary.  
Let $\epsilon, \eta > 0$ with $\epsilon \ll \eta$.  
Let
    $$ h(r) = \frac{\ln(r/ \d )}{\ln( \eta / \d )}. $$
    Then there exists $\d > 0$ sufficiently small and a smooth function $\ve \colon (0,\infty) \to \R$ which satisfies the following properties:
\begin{enumerate}
\item $\ve(r) = h(r)$ for all $r \leq \eta - \epsilon$.  
\item $\ve(r) = 1$ for all $r \geq \eta$.  
\item $0 \leq \ve'(r) \leq \e$ for all $r \in (\eta - \epsilon, \eta)$.  
\item $| \ve''(r) | \leq \e$ for all $r \in (\eta - \epsilon, \eta)$.
\item $r \ve''(r) + 2 \ve'(r) \geq - \e$ for all $r \in (\eta - \epsilon, \eta)$.
\end{enumerate}
\end{lemma}

\begin{lemma}[Smoothing at $r=\d$]\label{lem:smoothing at r=delta}
Suppose $\epsilon = C \d$ for some fixed $0 < C < 1/4$. Let
    $$ h(r) = \frac{\ln(r/ \d )}{\ln( \eta / \d )}. $$
Then there exist smooth functions $\fe, \ve \colon (0,\infty ) \to \R$ and constants $C_1$ and $C_2$ which depend only on $C$ which satisfy
\begin{enumerate}
\item $\fe(r) = \d$ for $r \leq \d - \epsilon$.
\item $\fe(r) = r$ for $r \geq \d + \epsilon$
\item $0 \leq \fe'(r) \leq 1$ for all $r$.
\item $0\leq \fe''(r) \leq \frac{1}{C\delta}$ for all $r$.
\item $\ve(r) = 0$ for $r \leq \d - \epsilon$.
\item $\ve(r) = h(r)$ for $r \geq \d + \epsilon$.
\item $|\fe(r) \ve'(r) | < C_1 \left( \frac{1}{\ln(\eta / \d)} \right)$ for all $r$.
\item $|\fe^2(r) \ve''(r)| < C_2 \left( \frac{1}{\ln(\eta / \d)} \right)$ for all $r$.
\end{enumerate}

\end{lemma}

\begin{proof}[Proof of Lemma \ref{lem:smoothing at r=delta}]
We first construct the function $\fe$ and verify conditions $(1)-(4)$, we then construct the function $\ve$ and verify conditions $(5)-(6)$, then we conclude by verifying inequalities $(7)-(8)$.  

\noindent {\bf Construction of $\fe$:}  We begin by making the substitution $t = \dfrac{r - (\d - \epsilon)}{2 \epsilon}$ and defining
\begin{equation*}
\psi(t)=\begin{cases}
e^{-1/t} & t > 0 \\
0 & t\le 0
\end{cases},
\end{equation*}
and
\begin{equation*}
\rho(t) = \frac{\psi(t)}{\psi(t) + \psi(1-t)}.
\end{equation*}
It is easy to check that $\psi$ and $\rho$ are both smooth and that $\rho$ satisfies the following properties
\begin{equation*}
\rho(t) = 0 \hskip 10pt (t \leq 0) \qquad \rho(t) = 1 \hskip 10pt (t \geq 1) \qquad 0 < \rho'(t) \hskip 10pt (0 < t < 1).
\end{equation*}
It follows immediately from the definition of $\rho$ that
\begin{equation*}
\rho(1-t) = 1 - \rho(t)
\end{equation*}
which implies that
\begin{equation}\label{eqn:integral of rho}
\int_0^1 \rho(t) \, dt = \frac{1}{2}.
\end{equation}
We then define $\fe(r)$ as follows:
\begin{equation*}
\fe(r)=\begin{cases}
\d & r \leq \d - \epsilon \\
\d + 2 \epsilon \int_0^t \rho(s) \, ds & \d - \epsilon < r < \d + \epsilon \\
r & r \geq \d + \epsilon
\end{cases},
\end{equation*}
where we recall that $t = \dfrac{r - (\d - \epsilon)}{2 \epsilon}$.  

It is clear from the construction that $\fe$ is continuous at $r = \d - \epsilon$, and integral \eqref{eqn:integral of rho} shows that $\fe$ is continuous at $r = \d + \epsilon$. 
Using the chain rule, one has
\begin{equation*}
\fe'(r) = \rho(t) \qquad \fe''(r) = \frac{1}{2 \epsilon} \rho'(t) \geq 0.
\end{equation*}
Since the function $\rho'(t)$ is maximized at $\rho'(1/2)=2$, this also gives us the upper bound
\begin{equation*}
\fe''(r) = \frac{1}{2 \epsilon} \rho'(t) \leq \frac{1}{\epsilon}=\frac{1}{C\d}.
\end{equation*}
We have thus verified conditions $(1)-(4)$.  
A direct calculation shows that all derivatives of $\rho$ are equal to $0$ at $t=0, 1$, which verifies that $\fe$ is smooth.

\vskip 10pt

\noindent {\bf Construction of $\ve$:} Let $\xi$ be a small positive parameter, and define the functions $\psi_\xi$ and $\rho_\xi$ by
\begin{equation*}
\psi_\xi(t)=\begin{cases}
e^{-1/(t-\xi)} & t > \xi \\
0 & t\le \xi
\end{cases},
\end{equation*}
and
\begin{equation*}
\rho_\xi(t) = \frac{\psi_\xi(t)}{\psi_\xi(t) + \psi_\xi(1- \xi - t)}
\end{equation*}
where $t$ is the same as for the definition of $\fe$.  
If one sets $\xi = 0$ they obtain the functions $\psi$ and $\rho$ from above.  
The reason why we introduce the parameter $\xi$ is because now the non-negative function $\rho_\xi$ will be equal to $0$ on the $\xi$-neighborhood of $t=0$, and will equal $1$ on the $\xi$-neighborhood of $t=1$.  

Define the function $q(t)$ by
\begin{equation}\label{eqn:definition of q}
q(t) = \rho_\xi(t) \frac{\d + \epsilon}{(\d - \epsilon) + 2 \epsilon t} + D \rho_\xi(t) (1- \rho_\xi(2t-1)).
\end{equation}
where $D$ is a constant which is chosen so that
\begin{equation}\label{eqn:integral of q}
\int_0^1 q(t) \, dt = \frac{h(\d + \epsilon)}{2 \epsilon h'(\d + \epsilon)}.
\end{equation}
Distributing this integral over \eqref{eqn:definition of q} and solving gives an easy way to obtain an explicit formula for $D$. 
\begin{align*}
    D=\frac{\frac{h(\d+\epsilon)}{2\epsilon h'(\d+\epsilon)}-\int_0^1 \rho_\xi(t)\frac{\d+\epsilon}{(\d-\epsilon)+2\epsilon t}dt}{\int_0^1\rho_\xi(t)(1-\rho_\xi(2t-1))dt}=\frac{\frac{(1+C)\ln(1+C)}{2C}-\int_0^1 \rho_\xi(t)\frac{1+C}{1-C+2Ct}dt}{\int_0^1\rho_\xi(t)(1-\rho_\xi(2t-1))dt}.
\end{align*}
When $\xi=0$, we have $\rho_\xi=\rho$ and the numerator of $D$ will be 
\begin{align*}
    &\frac{(1+C)\ln(1+C)}{2C}-\int_0^1 \rho(t)\frac{1+C}{1-C+2Ct}dt\\
    &=\int_\frac{1}{2}^1\frac{1+C}{1-C+2Ct}dt-\int_0^1 \rho(t)\frac{1+C}{1-C+2Ct}dt\\
    &=-\int_0^\frac{1}{2}\rho(t)\Bigl(\frac{1+C}{1-C+2Ct}-\frac{1+C}{1+C-2Ct}\Bigr) dt.
\end{align*}
Moreover, we can easily check that 
\begin{align*}
    0<\frac{1+C}{1-C+2Ct}-\frac{1+C}{1+C-2Ct}<1,
\end{align*}
and thus $-1<D<0$. By continuity, the same inequalities hold for all sufficiently small $\xi>0$.

We can now define $\ve$ by
\begin{equation*}
\ve(r)=\begin{cases}
0 & r \leq \d - \epsilon \\
2 \epsilon h'(\d + \epsilon) \int_0^t q(s) \, ds & \d - \epsilon < r < \d + \epsilon \\
h(r) & r \geq \d + \epsilon.
\end{cases},
\end{equation*}
By \eqref{eqn:integral of q} and the definition of $t$ it is immediate that $\ve$ is continuous at $r = \d - \epsilon$ and $r = \d + \epsilon$.  
Moreover, by introducing the parameter $\xi$, one can easily show that $\ve$ is smooth near $r = \d - \epsilon$ and $r = \d + \epsilon$ as follows.  

The function $\rho_\xi(t)$ is identically equal to $0$ over the $\xi$-neighborhood of $0$.  
For $\xi$ chosen sufficiently small, this corresponds to a small neighborhood of $r = \d - \epsilon$.  
A direct observation shows that $q = 0$ whenever $\rho_\xi = 0$, and thus $\ve$ is smooth near $r = \d - \epsilon$.  

Restricting $t$ near $1$ corresponds to $r$ near $\d + \epsilon$.  
For $t$ sufficiently close to $1$, $\rho_\xi(t) = 1$ and $\rho_\xi(2t-1) = 1$.  
Thus, for $t$ close to $1$, we have
\begin{equation*}
q(t) = \frac{\d + \epsilon}{(\d - \epsilon) + 2 \epsilon t}.
\end{equation*}
Applying the chain rule, we have
\begin{equation*}
\ve'(r) = h'(\d + \epsilon) q(t) = \frac{1}{(\d + \epsilon) \ln \left( \frac{\eta}{\d} \right)} \cdot \frac{\d + \epsilon}{r} = \frac{1}{r \ln \left( \frac{\eta}{\d} \right)} = h'(r).
\end{equation*}
Hence $\ve = h$ for $r$ sufficiently close to $\d + \epsilon$, and therefore $\ve$ is smooth.  
This verifies properties $(5) - (6)$. 

As an intermediate step, we need to show that $0 \leq \ve'(r) \leq C_1 \cdot h'(\d + \epsilon)$ for all $r$.  
The lower bound is immediate:
\begin{equation*}
\ve'(r) = h'(\d + \epsilon) q(t)  \geq 0
\end{equation*}
since $q(t)\geq \rho_\xi(t)(1+D)\geq 0$.

For the upper bound, recall that $\epsilon = C \d$ for some constant $0 < C \leq 1/4$.  
Because of this, each of the terms $\rho_\xi$, $D$, and $(\d + \epsilon)/(\d - \epsilon + 2 \epsilon t)$ are uniformly bounded above on the interval $r\in (\d - \epsilon, \d+\epsilon)$ (or $t\in (0,1)$) by constants which depend only on $C$.  
Therefore, there exists a constant $C_1$ depending only on $C$ such that $q(t) < C_1$ for $t \in [0,1]$.  
Thus, over this same region
\begin{equation*}
\ve'(r) = h'(\d + \epsilon) q(t) < C_1 h'(\d + \epsilon).
\end{equation*}

For property $(7)$ note that, over the interval $(\d - \epsilon, \d + \epsilon)$ for $r$:
\begin{equation*}
|\fe(r) \ve'(r)| \leq (\d + \epsilon) \cdot C_1 h'(\d + \epsilon) = \frac{C_1 (\d + \epsilon)}{(\d + \epsilon) \ln \left( \eta / \d \right)} = C_1 \left( \frac{1}{\ln(\eta / \d) } \right).
\end{equation*}

We conclude by verifying property $(8)$.  
For $r \in (\d - \epsilon, \d + \epsilon)$, we have
\begin{equation*}
\ve''(r) = \frac{h'(\d + \epsilon)}{2 \epsilon} q'(t) = \frac{1}{2 \epsilon (\d + \epsilon) \ln( \eta / \d)} q'(t).
\end{equation*}
Using the same reasoning as for $q$, we know that there exists a constant $K(C)$ depending only on $C$ such that $|q'(t)| < K$.
We then have
\begin{equation*}
|\fe^2(r) \ve''(r)| < \frac{K (\d + \epsilon)^2}{2 \epsilon (\d + \epsilon) \ln ( \eta / \d )} = \frac{K \d (1+C)}{2 C \d \ln(\eta / \d)} = C_2 \left( \frac{1}{\ln ( \eta / \d )} \right).
\end{equation*}

\end{proof}

\bibliographystyle{abbrv}
\bibliography{Bibliography}

\end{document}